\documentclass[a4paper,10pt]{article}

\usepackage{amsthm}
\usepackage{amsmath}
\usepackage{amsxtra}
\usepackage{amssymb}
\usepackage{thmtools}
\usepackage[colorlinks=true,linkcolor=blue,citecolor=blue,urlcolor=blue]{
hyperref}
\usepackage{mathrsfs}
\usepackage{undertilde}
\usepackage{turnstile}
\usepackage{enumitem}
\usepackage[retainorgcmds]{IEEEtrantools}
\declaretheorem[name=Definition,style=definition,qed=$\dashv$,
numberwithin=section]{dfn}

\declaretheorem[name=Theorem,style=plain,sibling=dfn]{tm}
\declaretheorem[name=Theorem,style=plain,numbered=no]{tm*}
\declaretheorem[name=Fact,style=plain,sibling=dfn]{fact}
\declaretheorem[name=Fact,style=plain,numbered=no]{fact*}
\declaretheorem[name=Lemma,style=plain,sibling=dfn]{lem}

\declaretheorem[name=Corollary,style=plain,sibling=dfn]{cor}
\declaretheorem[name=Remark,style=definition,sibling=dfn]{rem}
\declaretheorem[name=Claim,style=plain]{clm}
\declaretheorem[name=Claim,style=plain,numbered=no]{clm*}

\declaretheorem[name=Sublaim,style=plain,numbered=no]{sclm*}

\declaretheorem[name=Question,style=definition,sibling=dfn]{ques}
\declaretheorem[name=Conjecture,style=definition,sibling=dfn]{conj}

\newcommand{\Con}{\mathrm{Con}}

\newcommand{\Lim}{\mathrm{Lim}}

\newcommand{\QQ}{\mathbb Q}

\newcommand{\PP}{\mathbb P}

\newcommand{\sub}{\subseteq}
\newcommand{\cross}{\times}
\newcommand{\all}{\forall}

\newcommand{\inter}{\cap}
\renewcommand{\int}{\inter}

\newcommand{\om}{\omega}
\newcommand{\pow}{\mathcal{P}}
\newcommand{\OR}{\mathrm{OR}}

\newcommand{\Hull}{\mathrm{Hull}}

\newcommand{\cut}{\backslash}

\newcommand{\Ll}{\mathcal{L}}

\newcommand{\rg}{\mathrm{rg}}
\newcommand{\dom}{\mathrm{dom}}

\newcommand{\crit}{\mathrm{cr}}

\newcommand{\rest}{\!\upharpoonright\!}
\newcommand{\com}{\circ}

\newcommand{\sats}{\models}
\newcommand{\elem}{\preccurlyeq}

\newcommand{\AC}{\mathsf{AC}}
\newcommand{\DC}{\mathsf{DC}}

\newcommand{\HOD}{\mathrm{HOD}}
\newcommand{\HC}{\mathrm{HC}}
\newcommand{\ZFC}{\mathsf{ZFC}}
\newcommand{\ZF}{\mathsf{ZF}}

\newcommand{\Coll}{\mathrm{Col}}

\newcommand{\her}{\mathcal{H}}

\newcommand{\id}{\mathrm{id}}

\newcommand{\forces}{\dststile{}{}}

\newcommand{\trancl}{\mathrm{trancl}}
\newcommand{\Vop}{\mathrm{Vop}}

\newcommand{\ZFR}{\mathrm{ZFR}}

\newcommand{\OD}{\mathrm{OD}}

\newcommand{\psub}{\subsetneq}

\newcommand{\tu}{\textup}

\newcommand{\eot}{\mathrm{eot}}

\DeclareMathOperator{\Th}{Th}

\DeclareMathOperator{\card}{card}
\DeclareMathOperator{\cof}{cof}

\begin{document}

\title{Reinhardt cardinals and iterates of $V$}

\author{Farmer Schlutzenberg\footnote{schlutze@uni-muenster.de}}

\maketitle

\begin{abstract} Assume $\ZF(j)$ and there is a Reinhardt cardinal,
as witnessed by the elementary embedding $j:V\to V$.
We investigate the linear iterates $(N_\alpha,j_\alpha)$ of $(V,j)$,
and their relationship to $(V,j)$, forcing and definability,
including that for each infinite $\alpha$, every set is set-generic over 
$N_\alpha$, but $N_\alpha$ is not a set-ground.

Assume second order $\ZF$. We prove
that the existence of super Reinhardt cardinals and total Reinhardt
cardinals is not affected by small forcing.
And if $V[G]$ has a set of ordinals which is not in $V$,
then $V[G]$ has no 
elementary embedding $j:V[G]\to M\sub V$ (even allowing $M$
to be illfounded).\footnote{Keywords: Reinhardt cardinal,
iterate, HOD, definability, forcing, Axiom of Choice}
\footnote{\label{ftn:results_moved}
This paper originates from a set of rather informal notes 
\cite{reinhardt_non-definability}, published on 
arxiv.org (v1 of arXiv:2002.01215). The current paper is v$n$ of 
(the same) arXiv:2002.01215, where $n>1$. Some of the 
results from the original notes,
regarding definability and 
constructibility of embeddings $j:V_\delta\to V_\delta$ and related facts,
and those on extenders under ZF and definability of $V$-critical points from a 
proper class of weak L\"owenheim Skolem cardinals, are to appear
in the papers \cite{cumulative_periodicity} and \cite{ZF_extenders}.
The original notes together with some errata can be seen at
https://sites.google.com/site/schlutzenberg/home-1,
but it is in general better to refer to the current paper
and the two just cited. Of course, the original notes  also remain available
on arxiv.org, as mentioned above.}
\end{abstract}

\section{Introduction}

A \emph{Reinhardt cardinal}, introduced by William Reinhardt
in \cite{reinhardt_diss} and \cite{reinhardt_remarks},
is the critical point of an elementary embedding 
$j:V\to V$.
Kunen showed in \cite{kunen_no_R} that if $V\sats\ZFC$ then there is no such 
$j$.
So Reinhardt embeddings (that is, $j:V\to V$) are considered with background 
theory
$\ZF(j)$ (for the basic definitions \S\ref{subsec:notation} below), 
possibly augmented with fragments of
$\AC$. 
 We work throughout in $\ZF$
or variants thereof such as $\ZF(j)$ or $\ZF_2$, indicating any choice 
assumptions where they are adopted.

This paper primarily investigates the iterates of $V$ 
under an elemenetary $j:V\to V$.
Let $M_0=(V,j)$. In Definition \ref{dfn:iterates}
we define the $\alpha$th iterate
$M_\alpha=(N_\alpha,j_\alpha)$ of $M_0$ for each $\alpha\in\OR$.
These are analogous
to the iterates of $(L(V_{\lambda+1}),k)$ in the context of $I_0$,
investigated by Woodin in \cite{woodin_sem2}, and Mohammad Golshani 
posed a  question
of the nature of the $M_\alpha$ on the discussion board Mathoverflow.

In \cite{gen_kunen_incon}, many variants of and facts related to Kunen's 
$\ZFC$ inconsistency result are discussed. This paper makes a contribution 
toward understanding those issues under $\ZF$, a question posed in 
\cite{gen_kunen_incon}. The analysis of the iterates
$M_\alpha$ will show that 
some of the conclusions proven in \cite{gen_kunen_incon} fail in a 
strong way under $\ZF(j)+j:V\to V$.
On the other hand, we are able to establish a partial $\ZF$ generalization of 
Woodin's $\ZFC$ theorem \cite[Theorem 5]{gen_kunen_incon} excluding 
elementary $j:V[G]\to V$. 

We prove the following facts. In Theorem \ref{tm:every_set_generic_over_N_tau},
we show that every set is set-generic over $N_\alpha$.
The key to this is using a slight variant of a forcing due to Schindler
(the latter being related to the Bukowsky forcing, and also to Vopenka and 
the extender algebra).
In particular, for the iteration map $i_{0\alpha}:V\to N_\alpha$,
$i_{0\alpha}\rest V_\beta$ is set-generic over $N_\alpha$,
and this has the consequence that there are interesting
generic extensions of $V$, including that $\crit(j)$
is a \emph{virtual Berkeley cardinal} (see Definition 
\ref{dfn:virtually_Berkeley} and Corollary \ref{cor:forcing_extensions}).
In Theorem \ref{tm:N_tau_not_set-ground}, we show that, however, for 
$\alpha\geq\om$, $N_\alpha$ is not a set-ground of $V$ (note $N_n=V$ for 
$n<\om$).
It is well known that every set is also set-generic over $\HOD$.
We show in Theorem \ref{tm:lambda_Mahlo} that $\HOD\sub N_\om$, 
and in fact $\Hull^{(V,M_\om)}(N_\om)=N_\om$ (where the hull is computed
in $V$, with the predicate $M_\om$ (including $j_\om$) from parameters in 
$N_\om$), but there are $\alpha\in\OR$ with $\HOD\not\sub 
N_\alpha$.

By \cite[Theorem 12]{gen_kunen_incon},
 under $\ZFC_2$, if $j:V\to M$ is elementary with $M$ transitive then 
$V=\Hull^V(\rg(j))$.
This fails trivially in the case that $(V,j)\sats\ZFR$ (where $M=V$),
since then $\Hull^V(\rg(j))=\rg(j)$. However, because $\Hull^V(N_\om)=N_\om$,
we get a non-trivial failure, as we have the iteration map $i_{0\om}:V\to 
N_\om$, and $N_\om\neq V$. We will also prove a stronger fact in Theorem 
\ref{tm:V=HOD(X)_no_R}, that for every set $X$, 
$V\neq\Hull^{(V,M_\om)}(N_\om\cup X)$.
At the same time we will establish
a strengthening of Suzuki's theorem  \cite{suzuki_no_def_j} 
that $j$ is not definable from parameters (which was strengthened
in other ways in \cite{gen_kunen_incon}, \cite{cumulative_periodicity} and 
\cite{ZF_extenders}):
we show that $j\rest\OR$ is not definable from parameters and the predicate 
$(N_\om,j_\om)$.\footnote{\label{ftn:j_rest_OR_not_def}The fact that 
$j\rest\OR$ is not definable 
over $V$ from parameters (but without any extra predicate)
is already a consequence of \cite[Theorem 35]{gen_kunen_incon},
by which  no elementary $k:\HOD\to\HOD$ is definable 
from parameters (without assuming $\AC$).
In Remark \ref{rem:j_rest_OR_not_def_from_params}
we also sketch an alternate proof of the latter fact.}

It follows in particular that if there is a Reinhardt,
then for all sets $X$,
 $V\neq\HOD(X)$.
 And if $V$ is total Reinhardt or there is a Berkeley cardinal,
a similar proof gives that
$V\neq\HOD_A(X)$ for any class $A$ and set $X$.
But in fact much more is now known;
see Remark \ref{rem:goldberg_usuba}.

In \S\ref{sec:small_forcing} we prove some facts regarding small 
forcing and (very) large cardinals.
Theorem 
\ref{tm:Reinhardt_with_enough_AC_restricts_to_ground} asserts that if $\PP\in 
V$ and $G$ is $(V,\PP)$-generic
and $(V[G],j)\sats\ZFR+$``$\PP$ is wellorderable and $\card(\PP)^+$-$\DC$ 
holds'',
then $(V,j\rest V)\sats\ZFR$.
Assume $\ZF_2$. By Theorem \ref{tm:no_new_sR_or_tR}, if $G$ is 
set-generic
over $V$ and $V[G]$ has a super Reinhardt cardinal,
then so does $V$. A similar fact is proven for total Reinhardtness.
 Theorem \ref{tm:enough_AC_ground_def} is a variant
of the $\ZF+\DC_\delta$ ground definability result of Gitman-Johnstone 
\cite{ground_def_gitman_johnstone},
in which we assume some $\DC$ in the generic extension,
but not  (at least not explicitly)  in the ground model
 (but we assume that the forcing
in question is small).

Theorem \ref{tm:no_gen_emb_into_V_with_new_O_set} is a generalization
of Woodin's $\ZFC_2$ theorem \cite[Theorem 4]{gen_kunen_incon} that if $G$ 
is set-generic over $V$
then there is no elementary $j:V[G]\to V$. We prove
 assuming $\ZF_2$ that if $G$ is set-generic
 over $V$ \emph{and $V[G]$ has some set of ordinals 
not in $V$}, then there is no elementary $j:V[G]\to M\sub V$,
even allowing $M$ to be illfounded.
This fact appears to be new even in the $\ZFC$ 
context; the results in \cite{gen_kunen_incon}
seem to rely on stationary set combinatorics which does
not seem to be available in this generality.

There are many questions that are open,
and some of these are mentioned throughout.
We would like to thank Gabriel Goldberg and Toshimichi
Usuba for their feedback regarding certain questions
in the notes \cite{reinhardt_non-definability}
(upon which the material here is based), and other helpful suggestions.

\subsection{Basic definitions}\label{subsec:notation}

$\OR$ denotes the class of ordinals and $\Lim$ the limit ordinals.

The \emph{language $\Ll_{A_1,\ldots,A_n}$ of set theory with predicates 
$A_1,\ldots,A_n$}
is the first order language with binary 
predicate symbol
$\in$ and predicate symbols $A_1,\ldots,A_n$.
The theory \emph{$\ZF(\vec{A})$} is the theory in 
$\Ll_{\in,\vec{A}}$
with all $\ZF$ axioms, allowing all 
formulas of $\Ll_{\in,\vec{A}}$
in the Separation and Collection schemes (so each $A_i$ represents a class).

$\ZFR$ ($\ZF+$ Reinhardt) is the theory 
$\ZF(\widetilde{j})+$``$\widetilde{j}:V\to V$ is 
$\Sigma_1$-elementary''.
And $\ZFR(A)$ is the theory 
$\ZF(\widetilde{j},A)+$``$\widetilde{j}:(V,A)\to(V,A)$ is 
$\Sigma_1$-elementary''.
Essentially by \cite[Proposition 5.1]{kanamori}, $\ZFR(A)$ proves
$\widetilde{j}$ is fully elementary (as a theorem scheme).

Second order set theory is denoted $\ZF_2$ (see 
\cite{woodin_koellner_bagaria}).
  Models are of the form $(V,\in,P)$,
  where $(V,\in)\sats\ZF$ and $P$ is a collection of classes/subsets of 
$V$,
  satisfying the $\ZF_2$ axioms.
Given a transitive $W$ and $P\sub\pow(W)$,
 we write $(W,P)\sats\ZF_2$ iff $(W,\in,P)\sats\ZF_2$.
When we ``work in $\ZF_2$'', we mean that we work in
such a model $W$, and all talk of proper classes refers to elements of $P$.

 Let $\delta\in\Lim$.
The \emph{cofinality} of $\delta$, \emph{regularity}, \emph{singularity}
are defined as usual (in terms of cofinal functions between ordinals).
 We say  $\delta$ (or $V_\delta$) is \emph{inaccessible}
 iff there is no $(\gamma,f)$ such that $\gamma<\delta$
 and $f:V_\gamma\to\delta$ is cofinal in $\delta$.

Work in $\ZF_2$. 
Then
 $\kappa\in\OR$ is \emph{Reinhardt}
iff there is a class $j$ such that $(V,j)\sats\ZFR$ and $\kappa=\crit(j)$.
Following \cite{woodin_koellner_bagaria},
 $\kappa\in\OR$ is \emph{super-Reinhardt}
 iff for every $\lambda\in\OR$ there is a class $j$
 such that $(V,j)\sats\ZFR$ and $\crit(j)=\kappa$
 and $j(\kappa)\geq\lambda$. And $\OR$ is \emph{total Reinhardt}
 if for every class $A$ there is $\kappa\in\OR$ such that
$\kappa$ is \emph{$({<\OR},V,A)$-reflecting},
meaning that for every $\lambda\in\OR$ there is a 
class $j$ such that 
$(V,j,A)\sats\ZFR(A)$ and $\crit(j)=\kappa$ and $j(\kappa)\geq\lambda$.
And $\delta$ is \emph{Berkeley} if for all transitive sets $M$
 with $\delta\in M$, and all $\eta<\delta$,
 there is an elementary $j:M\to M$ 
with $\eta<\crit(j)<\delta$.

Let $j:V_\delta\to V_\delta$ or $j:V\to V$ be $\Sigma_1$-elementary,
for a limit $\delta$ or $\delta=\OR$.
We write $\kappa_0(j)=\crit(j)$ and 
$\kappa_{n+1}(j)=j(\kappa_n(j))$
and $\kappa_\om(j)=\sup_{n<\om}\kappa_n$.

Let $C$ be a class. Then $\OD_C(X)$
denotes the class of all sets $y$ such that $y$
is definable from the predicate $C$
and parameters in $\OR\cup X$.
And $\HOD_C(X)$ denotes the class of all $y$
such that the transitive closure of $\{y\}$
is $\sub\OD_C(X)$.
And $\OD(X)=\OD_\emptyset(X)$ and $\HOD(X)=\HOD_\emptyset(X)$.
And for a structure $M$ and $X\sub M$, $\Hull^M(X)$ denotes
the collection of elements $x\in M$ such that $x$ is definable
over $M$ from (finitely many) parameters in $X$.

\section{Generic $j:V[G]\to M\sub V$}\label{sec:generic_j:V[G]_to_M}

The following fact, due to Woodin, is proved in \cite[Theorem 
5]{gen_kunen_incon}:
\begin{fact}[Woodin]\label{fact:Woodin_j:V[G]_to_V} Assume $\ZFC$. Let $G$ be 
set-generic over $V$.
Then there is no elementary $j:V[G]\to V$.
\end{fact}
Of course in the case that $G\in V$, this reduces
to Kunen's inconsistency.
Various generalizations of this result can also be seen in 
\cite{gen_kunen_incon}; one in particular is that,
by \cite[Corollary 34]{gen_kunen_incon}, under $\ZF$, if $G$ is set-generic 
over $V$,
then there is no elementary $j:V[G]\to V$ which is definable from parameters 
over $V[G]$.

In this section, we generalize Fact \ref{fact:Woodin_j:V[G]_to_V} 
in a few ways. First, we replace
the assumption that $\ZFC$ holds in $V$ with $\ZF$
plus 
the requirement that $G$ adds a set of ordinals to $V$
(which of course holds in the $\ZFC$ context, if $G\notin V$).
Second, we only require that $M\sub V$, not that $M=V$
or $M$ is eventually stationarily correct or anything further,
and so in particular, we allow  $M$ to be illfounded. 
The proof also gives an alternate proof of
Woodin's result, for the case that $G\notin V$;
note that we do not use any stationary set combinatorics.
The same kind of argument can be used in certain other places
where the stationary set argument is traditionally used, for example,
in Usuba's $\ZFC$ proof that the mantle coincides with the $\kappa$-mantle if 
$\kappa$ is extendible.

\begin{lem}\label{lem:no_gen_emb_into_V_with_new_O_set} Assume $\ZF$.
 Let $G$ be $(V,\PP)$-generic for some forcing $\PP\in V$
and let $\alpha\in\OR$.
Let $\theta\in\OR$ be such that $V$ has no surjection
$\PP^{<\om}\cross\alpha\to\theta$.
Suppose $(M,E,N,j)$ are such that
$M,E\in V$ and $N,j\in V[G]$, $N$ is transitive and rudimentarily closed,
$\her_\theta^V\sub N$, and
$j:(N,\in)\to(M,E)$
is elementary.
\tu{(}Note we allow $(M,E)$ to be illfounded.\tu{)} Then:
\begin{enumerate}
\item\label{item:j_rest_alpha_in_V} $j\rest\alpha\in V$,
\item\label{item:pow(alpha)_sub_V} $\pow(\alpha)\inter N=\pow(\alpha)\inter V$.
\end{enumerate}
\end{lem}
\begin{proof}
Part \ref{item:pow(alpha)_sub_V}: This follows directly from part 
\ref{item:j_rest_alpha_in_V}. For $\pow(\alpha)\inter V\sub\her_\theta^V\sub 
N$ by assumption, and supposing $X\in N\inter\pow(\alpha)$,
 for each $\xi<\alpha$, we have
\[\xi\in X\iff N\sats\text{``}\xi\in X\text{''}\iff
(M,E)\sats\text{``}j(\xi)\in j(X)\text{''},\]
but since $(M,E)\in V$ and $j\rest\alpha\in V$ (by part 
\ref{item:j_rest_alpha_in_V}), therefore $X\in V$.

Part \ref{item:j_rest_alpha_in_V}: We prove a couple of claims.

 \begin{clm}\label{clm:no_surjection_to_theta} $V[G]\sats$``there is  no 
surjection $f:\PP^{<\om}\cross\alpha\to\theta$''.\end{clm}
\begin{proof} Suppose otherwise,
 and fix $\widetilde{f}\in V$  with $\widetilde{f}_G=f$
 and $p_0\in G$ such that
$p_0\forces_\PP$``$\widetilde{f}:\check{\PP}^{<\om}\cross\check{\alpha} 
\to\check { \theta}$ is surjective''.
We define a function $g\in V$. For
$(p,\vec{q},\xi,\gamma)$ in $\PP\cross\PP^{<\om}\cross\alpha\cross\theta$, set
 \[ g(p,\vec{q},\xi)=\gamma\iff
p\leq 
p_0\text{ 
and }p\forces_\PP\widetilde{f}(\check{\vec{q}},\check{\xi})=\check{\gamma}. \]
Note $A=\dom(g)\sub\PP\cross\PP^{<\om}\cross\alpha$
and $g:A\to\theta$ is surjective, impossible.
\end{proof}

In $V[G]$, we have
$j\rest\theta:\theta\to\OR^M$.
Fix $\widetilde{k}\in V$ with $\widetilde{k}_G=j\rest\theta$ 
and $p_0\in G$ with
$p_0\forces_\PP$``$\widetilde{k}:\check{\theta}\to\OR^{\check{M}}$''.
\footnote{Here $\OR^M=\{x\in M\bigm|(M,E)\sats\text{``}x\text{ is an 
ordinal''}\}$, and
 $M$ might be illfounded.}
Working in $V$, for  $p\in\PP$ with $p\leq p_0$, define
$X_p\sub\theta$ and $j_p:X_p\to\OR^M$ the function giving all values of
$\widetilde{k}$ decided by $p$; that is, for $\xi<\theta$,
\[ j_p(\xi)=x\iff x\in\OR^{M}\text{ and }p\forces_\PP
\widetilde{k}(\check{\xi})=\check{x},\]
and $X_p=\dom(j_p)$. So $X_p,j_p\in V$.

\begin{clm}\label{clm:X_p_long} There is $p\in G$ such that $p\leq p_0$ and 
$X_p$ has ordertype $\beta>\alpha$.\footnote{Note that the proof does not use 
$\AC$.}
\end{clm}
\begin{proof}We have
$\theta=\bigcup_{p\in G, p\leq p_0}X_p$.
Suppose the claim fails.
In $V[G]$,  define a surjection
$f:\PP\cross\alpha\to\theta$
by setting $f(p,\xi)=$ the $\xi$th element of $X_p$,
in the case that $X_p$ is defined and has ordertype $>\xi$,
and $f(p,\xi)=0$ otherwise.
This contradicts Claim \ref{clm:no_surjection_to_theta}.\end{proof}

So fix $p$ as in Claim \ref{clm:X_p_long}.
Let $\alpha'$ be the $\alpha$th element of $X_p$.
So $Y=X_p\inter\alpha'$ has ordertype $\alpha$.
Let $\pi:\alpha\to Y$ be the increasing enumeration of
$Y$.
Then $\alpha,Y,\pi\in\her_\theta^V\sub N$.
Since $N$ is transitive,
\[ N\sats\text{``}\pi:\alpha\to Y\text{ is the increasing enumeration of  }
Y\text{''},\]
so applying $j$,
\[ (M,E)\sats\text{``}j(\pi):j(\alpha)\to j(Y)\text{ is
 the increasing enumeration
of }j(Y)\text{''}.\]

Let $j_Y=j\rest Y\in V=j_p\rest Y\in V$. Note that
$j_Y(\pi(\xi))=j(\pi(\xi))=j(\pi)(j(\xi))$
for each $\xi<\alpha$,\footnote{Here $j(\pi)(j(\xi))$
is a slight abuse of notation; it should really be $(j(\pi)(j(\xi)))^M$,
i.e. the $x\in M$ such that $M\sats$``$x=j(\pi)(j(\xi))$''.} and so\footnote{
The notation ``$j(\pi)^{-1}$'' involves a similar abuse.}
$j\rest\alpha=j(\pi)^{-1}\com j_Y\com \pi$,
and since $j_Y,\pi,j(\pi),(M,E)\in V$, we get $j\rest\alpha\in 
V$.
\end{proof}

From the lemma, we can deduce:

\begin{tm}\label{tm:no_gen_emb_into_V_with_new_O_set}
Assume $(V,P)\sats\ZF_2$.
 Let $G$ be $(V,\PP)$-generic for some $\PP\in V$.
Suppose there are $\alpha\in\OR$ and $X\sub\alpha$ with $X\in V[G]\cut V$.

Then there is no $(j,M)\in P[G]$ with $M\sub V$
and $j:V[G]\to M$ elementary.
Hence, there is no $j\in P[G]$
such that $j:V[G]\to V$ is elementary.
\end{tm}

\begin{ques}
  Is it consistent that there is
 $j:V[G]\to M\sub V$
 as in the statement of Theorem \ref{tm:no_gen_emb_into_V_with_new_O_set}, if 
we 
instead have
 $\pow(\alpha)\inter V[G]\sub V$
 for each $\alpha\in\OR$? Does it follow that
 $V^{<\OR}\inter V[G]\sub V$?
What about in the case that $M=V$? 

What about the converse direction, i.e. an embedding $j:V\to V[G]$,
where $G$ is set-generic over $V\sats\ZF$?
\end{ques}

In the next section we will see that if $\ZFR$ is consistent,
then one \emph{can} get $j:V\to M\psub V$ where every \emph{set} is 
(individually) set-generic over $M$; this fact should be compared with 
\cite[\S4]{gen_kunen_incon}.

\section{The iterates of $(V,j)$}\label{sec:no_Reinhardt_in_HOD(X)}

\subsection{Basic properties}

\begin{rem} As the reader familiar with inner model theory
will see, the methods in this section are heavily based on
analyses of $\HOD$ via direct limit systems of mice,
and related methods from the study of Varsovian models.
The present context actually provides a nice introduction
to some of those ideas, which avoids various  issues
(iteration trees, fine structure, comparison).

We start 
with an easy observation.
 Assume $(V,j)\sats\ZFR$
 and let $\lambda=\kappa_\om(j)$.
Then for each $X\in V_\lambda$
(i) there is no cofinal $f:\om\to\lambda$ with $f\in\OD_X$,
and (ii) $j\rest\lambda\notin\OD_X$.

For suppose $f$ is $\OD_X$. Then there is in fact some such $f$
which is definable from $(X,\lambda)$, by minimizing on other ordinal 
parameters.
Taking $n<\om$ large enough, then $j^n(X)=X$, so $j^n(f)=f$,
so $j^n(f(k))=f(k)$ for each $k<\om$.
But $\rg(f)$ is unbounded in $\lambda$,
and $\lambda$ is the least fixed point of $j^n$
which is $>\crit(j^n)$, a contradiction.
And from $j\rest\lambda$ we can define $\left<\kappa_n\right>_{n<\om}$,
which we have just seen is not $\OD_X$, so $j\rest\lambda$
is not $\OD_X$.
\end{rem}

Now in the study of $L(V_{\lambda+1})$ under $I_0$,
(see \cite{woodin_sem2}, for example), the 
iterates of $(L(V_{\lambda+1}),j)$
are an important focus.
And Mohammad Golshani asked, on the discussion board 
Mathoverflow,
what properties might hold of the intersection of the models
$N_\alpha$ which we now define.\footnote{\label{ftn:mof}For 
the discussion see
https://mathoverflow.net/questions/185253/reinhardt-cardinals-and-iterability.
Note there is a slight difference in notation;
we are writing $N_\alpha$ for what 
he refers to as $M_\alpha$ there.}

\begin{dfn}[Iterates $M_\alpha$]\label{dfn:iterates}
Assume $(V,j)\sats\ZFR$.
Set $M_0=(V,j)$.
We define, for $\alpha\leq\beta\in\OR$, the $\alpha$th iterate 
$M_\alpha=(N_\alpha,j_\alpha)$, where $N_\alpha$
is a transitive proper class and $j_\alpha:N_\alpha\to N_\alpha$
a Reinhardt embedding of $N_\alpha$ (so $(N_\alpha,j_\alpha)\sats\ZFR$)
and $i_{\alpha\beta}:M_\alpha\to M_\beta$ is elementary
(literally a class function $N_\alpha\to N_\beta$, but elementary
with respect to the predicates $j_\alpha,j_\beta$), and such that
$i_{\alpha\gamma}=i_{\beta\gamma}\com i_{\alpha\beta}$
for all $\alpha\leq\beta\leq\gamma$, and $i_{\alpha\alpha}=\id$,
as follows. Given $M_\beta$, set
$M_{\beta+1}=(N_\beta,j_\beta(j_\beta))$
and $i_{\alpha,\beta+1}=j_\beta\com i_{\alpha\beta}$.
Given $M_\beta$ for all $\beta<\eta$ where $\eta$ is a limit,
set
$M_\eta$ as the direct limit,
and $i_{\alpha\eta}$ the direct limit map.

By elementarity, $N_\alpha$ is extensional.
Assuming that $N_\alpha$ is wellfounded, we identify it with its transitive 
collapse.

We introduce some symbols to the $\ZFR$ language to express these notions 
in formal language. For example, $\widetilde{M}$ represents the function 
$\alpha\mapsto M_\alpha$ defined under $\ZFR$ as above. So with the models
above and  $\alpha\in\OR^V$, we have 
$(\widetilde{M}_\alpha)^{(V,j)}=M_\alpha$,
and if $M_\alpha$ is wellfounded then
$(\widetilde{M}_\beta)^{M_\alpha}=M_{\alpha+\beta}$. We write 
$\widetilde{N}_\alpha,\widetilde{j}_\alpha,\widetilde{i}_{\alpha\beta}$
similarly. (This literally overloads the symbol 
$\widetilde{j}$ with two interpretations, but when we write $\widetilde{j}$
without a subscript, it always refers to the original function 
$\widetilde{j}_0$.) 
\end{dfn}

We now work in a model $(V,j)\sats\ZFR$.

As pointed out by Hamkins on Mathoverflow (see Footnote \ref{ftn:mof}),
the usual proof of linear iterability with respect to a single
normal measure gives the following (iterability
also holds in the $I_0$ case and generalizations thereof,
but the proof is finer there;
see \cite{woodin_sem2}
and \cite{con_lambda_plus_2}):
\begin{lem}
 For all $\alpha\in\OR$, $M_\alpha$ is a wellfounded
\tu{(}we take it transitive\tu{)},
$M_\alpha=(N_\alpha,j_\alpha)\sats\ZFR$, and $i_{\beta\alpha}:M_\alpha\to 
M_\beta$ is elementary
 and $\in$-cofinal.\footnote{Note that the elementarity
 is with respect to $M_\alpha=(N_\alpha,j_\alpha)$ and 
$M_\beta=(N_\beta,j_\beta)$, not just $N_\alpha,N_\beta$.}
\end{lem}
\begin{proof}
Everything is proved by induction. We get the elementarity
because $M_0=(N_0,j_0)\sats\ZFR$, so for each $n<\om$
there is a club class of $\alpha$ such that $(V_\alpha^{N_0},j\rest 
V_\alpha^{N_0})\elem_n(N_0,j)$, which easily gives that $i_{01}$
is elementary with respect to $M_0,M_1$, not just $N_0,N_1$, and 
likewise for larger indices.

Now suppose that $\eta\in\OR$ and $M_\eta$ is 
illfounded, and $\eta$ is least such. Then $\OR^{M_\eta}$ is illfounded 
(consider ranks of elements).
Since $M_{\alpha+1}$ and $M_\alpha$ have the same universe, $\eta$
is a limit ordinal. Let $\lambda\in\OR$ be least with $i_{0\eta}(\lambda)$ 
in the 
illfounded part of $\OR^{M_\eta}$.
There is $\alpha<\eta$ and $\lambda_1\in\OR$
such that $\lambda'=i_{0\alpha}(\lambda)>\lambda_1$
and $i_{\alpha\eta}(\lambda_1)$ is in the illfounded part of $M_\eta$.
Let $\eta'=i_{0\alpha}(\eta)\geq\eta$.
By  elementarity  of $i_{0\alpha}$,
$M_\alpha\sats$``$\eta'$ is the least 
$\gamma\in\OR$ such that $\widetilde{M}_\gamma$ is 
illfounded''.
Since $M_\alpha$ is wellfounded,
$(\widetilde{M}_\beta)^{M_\alpha}=M_{\alpha+\beta}$
and $\alpha+\eta'=\eta$, so $\eta'=\eta$ and 
$\widetilde{M}_\eta^{M_\alpha}=M_\eta$.
By elementarity, $M_\alpha\sats$``$\lambda'$
is the least $\xi\in\OR$ with $\widetilde{i}_{0\eta}(\xi)$
 in the illfounded part of $\widetilde{M}_\eta$''.
But $\widetilde{i}^{M_\alpha}_{0\eta}=i_{\alpha\eta}$,
and we chose $\alpha$ with $i_{\alpha\eta}(\lambda_1)$
in the illfounded part of $M_\eta$, a contradiction.
\end{proof}

\begin{lem}\label{lem:all_x_in_M_om_es}
Let $x\in M_\om$ and $m<\om$ with $x\in\rg(i_{m\om})$.
Then $i_{mn}(x)=x$ for all $n\in[m,\om)$.
\end{lem}
\begin{proof}
Let $\bar{x}$ be such that $i_{m\om}(\bar{x})=x$.
Then $M_m\sats$``$x=\widetilde{i}_{0\om}(\bar{x})$'',
so letting $(\bar{x}',x')=i_{mn}(\bar{x},x)$, elementarity gives 
$M_n\sats$``$x'=\widetilde{i}_{0\om}(\bar{x}')$'',
but $\widetilde{i}^{M_n}_{0\om}=i_{n\om}$, so by commutativity,
$x'=i_{m\om}(\bar{x})=x$.
\end{proof}

By the lemma, the following definition makes sense:

\begin{dfn} For $x\in M_\om$ let
$x^*=\lim_{n<\om}i_{n\om}(x)=i_{m\om}(x)$ for any/all $m<\om$
with $x\in\rg(i_{m\om})$.
\end{dfn}

\begin{lem}\label{lem:alpha*=}
For  $x\in M_\om$, we have
$x^*=i_{\om,\om+\om}(x)=\widetilde{i_{0\om}}^{M_\om}(x)$.
\end{lem}
\begin{proof}
Given $x$, let $\bar{x}$ and $m<\om$ with $i_{m\om}(\bar{x})=x$.
Then $M_m\sats$``$\widetilde{i}_{0\om}(\bar{x})=x$'',
so since $x^*=i_{m\om}(x)$, by elementarity we have
$M_\om\sats$``$\widetilde{i}_{0\om}(x)=x^*$'',
but $\widetilde{i}_{0\om}^{M_\om}=i_{\om,\om+\om}$,
so $i_{\om,\om+\om}(x)=x^*$.
\end{proof}

We now establish that $N_\om$ is closed under definability in $V$.
Note that the notation $\Hull^{(V,M_\om)}(X)$ refers to definability
in $V$ with $M_\om$ as a predicate, from parameters in $X$.

\begin{tm}\label{tm:lambda_Mahlo}
Let $\kappa_n=\crit(j_n)$ and $\lambda=\kappa_\om(j)$. Then:
\begin{enumerate}
 \item\label{item:HOD(Mom)} 
$\HOD\sub\Hull^V(N_\om)=\Hull^{(V,M_\om)}(N_\om)=N_\om$,
 \item\label{item:V_lambda^HOD} $V_\lambda^{\HOD}=V_{\lambda}^{\HOD^{N_\om}}$
 and $V_{\lambda+1}^{\HOD}\sub V_{\lambda+1}^{\HOD^{N_\om}}$.
 \item\label{item:Mahlo_in_HOD} $\lambda$ is greatly Mahlo in $\HOD$ and weakly 
compact in $\HOD^{N_\om}$.
\item\label{item:Nom-Berkeley} $\lambda$ is $N_\om$-Berkeley in 
$V$.\footnote{Note
this
generalizes
the (folklore?) fact that $\lambda$ is $\HOD$-Berkeley,
mentioned in \cite[\S ***]{cumulative_periodicity}.}
\item\label{item:lambda^+} 
$(\lambda^+)^{\HOD}\leq(\lambda^+)^{N_\om}<\lambda^+$
 \item If $\kappa_0^+$ is regular or $\cof(\kappa_0^+)=\kappa_0$
then $\cof((\lambda^+)^{M_\om})=\om$.
\item\label{item:if_cof<kappa_0} If $\cof(\kappa_0^+)<\kappa_0$ then 
$\cof((\lambda^+)^{M_\om})=\cof(\kappa_0^+)$.
\end{enumerate}
\end{tm}
\begin{proof}
Part \ref{item:HOD(Mom)}: We just need to see that 
$\Hull^{(V,M_\om)}(N_\om)\sub 
N_\om$.
We prove this by induction on rank of elements of 
$\Hull^{(V,M_\om)}(N_\om)$.
So let 
$X\in\Hull^{(V,M_\om)}(N_\om)$ with $X\sub N_\om$;
we must see that $X\in N_\om$. 
Let $\varphi$ be a formula in $\Ll_{\dot{M}}$ and $p\in N_\om$ be such that for 
$x\in N_\om$,
we have
$x \in X$ iff $V\sats\varphi(x,p;M_\om)$,
with $M_\om$ interpreting the predicate symbol $\dot{M}$.
Then for $n<\om$,
 $x\in X$ iff $M_n\sats\varphi(x,p;\widetilde{M}_\om)$
 iff $M_\om\sats\varphi(x^*,p^*;\widetilde{M}_\om)$
iff 
$M_\om\sats\varphi(\widetilde{i}_{0\om}(x),\widetilde{i}_{0\om}(p),\widetilde{M}
_\om)$, so $X$ is computed by $M_\om$, so $X\in N_\om$, as desired.

Part \ref{item:V_lambda^HOD}:
Because $\lim_{n<\om}\crit(i_{n\om})=\lambda$,
we easily have $V_\lambda^{\HOD}=V_{\lambda}^{\HOD^{N_\om}}$.
If $A\in V_{\lambda+1}^{\HOD}$ then just note that $A^*\in\HOD^{N_\om}$
and $A^*\inter V_\lambda=A$, so $A\in V_{\lambda+1}^{\HOD^{N_\om}}$.

Part \ref{item:Mahlo_in_HOD}: Since $\lambda=\crit(j_\om)$
and $j_\om\rest\HOD^{N_\om}:\HOD^{N_\om}\to\HOD^{N_\om}$ is elementary,
we easily get that $\lambda$ is weakly compact in $\HOD^{N_\om}$,
hence greatly Mahlo in $\HOD^{N_\om}$,
which by part \ref{item:V_lambda^HOD} implies $\lambda$
is greatly Mahlo in $\HOD$.

Part \ref{item:Nom-Berkeley}: 
By Lemma \ref{lem:all_x_in_M_om_es}, for every $M\in N_\om$, there is $n<\om$
such that $j_m(M)=M$ for all $m\in[n,\om)$,
and since $\lambda=\sup_{m<\om}\crit(j_m)$, this suffices.

Parts \ref{item:lambda^+}--\ref{item:if_cof<kappa_0}: These facts now
follow easily 
by considering the direct limit producing
$(\lambda^+)^{M_\om}$.
We have $(\lambda^+)^{M_\om}<\lambda^+$ because
there is a surjection from $\om\cross\lambda\to(\lambda^+)^{M_\om}$
computable directly from the direct limit.
If $\kappa_0^+$ is regular or $\cof(\kappa_0^+)=\kappa_0$
then  $j_{mn}$ (where $m\leq n<\om$) is discontinuous
at $\kappa_m^+$, and since $(\lambda^+)^{M_\om}$ is the direct limit
of an $\om$-sequence of these, this ordinal has cofinality $\om$.
In the other case, $j_{mn}$ and $j_{m\om}$ is continuous
at $(\kappa_0^+)$.
\end{proof}

\begin{ques}Is/can we have $V_{\lambda+1}^{\HOD}=V_{\lambda+1}^{\HOD^{N_\om}}$?
Or $V_{\lambda+1}^\HOD=V_{\lambda+1}^{N_\om}$? Is 
$\lambda$ Woodin in $\HOD$? Weakly compact in $\HOD$? 
Reflecting in $\HOD$?\footnote{Recall that $\kappa$ is a reflecting cardinal
iff $\kappa$ is inaccessible and for every stationary $A\sub\kappa$
there are stationarily many regular $\bar{\kappa}<\kappa$
such that $A\inter\bar{\kappa}$ is stationary in $\bar{\kappa}$.}
These properties seem (to the author)
reasonable to ask about, since $\lambda$ has them in $M_\om$
and $\HOD\sub M_\om$.
If $V_\lambda\sub\HOD$ then $\lambda$ is clearly
Woodin in $\HOD$, since it is Woodin (and more) in $N_\om$.
\end{ques}

\subsection{Every set is generic over $N_\alpha$}\label{sec:every_set_generic}

By Vopenka, every set is set-generic over $\HOD$. Moreover:

\begin{lem}\label{tm:every_om}
Every set in $V$ is set-generic
 over $M_\om$.
\end{lem}
\begin{proof}
One way to prove this is to use a Vopenka argument, as done in 
\cite{con_lambda_plus_2}. But we want to generalize the 
theorem
later to $M_\alpha$ for arbitrary $\alpha$, and there we will not have 
$\HOD\sub M_\alpha$. So we use a different argument, which seems to be
more generalizable.

We will show directly that for each $\beta\in\OR$
and each $X\sub V_\beta^{M_\om}$, with $X\in V$,
we can add $X$ set-generically to $M_\om$.
This suffices, because letting $\alpha\in\OR$, 
we can apply this to
 $X=i_{0\om}``V_\alpha$,
and then $V_\alpha^{M_\om[G]}=V_\alpha$
(and $G$ adds $j\rest V_\alpha$ also).

The forcing $\PP$ we use is the natural adaptation of Schindler's forcing
(to appear in \cite{vm2})
to the current setting (and Schindler's forcing is related to the forcings of 
Vopenka, Bukowsky and Woodin's extender algebra;
cf.~\ref{rem:Vop}).

Work in $M_\om$.
Let $\Ll$ be the class
of all infinitary Boolean formulas $\varphi$ in
propositional symbols $P_x$, for each $x\in V_\beta$.
That is, $\Ll$ is the minimal class generated by the following rules:
\begin{enumerate}
 \item For each $x\in V_\beta$, we have a corresponding propositional symbol 
$P_x\in\Ll$
 (with $x\neq y\implies P_x\neq P_y$).
 \item If $\varphi\in\Ll$ then $\neg\varphi\in\Ll$.
 \item If $A$ is a set and $A\sub\Ll$ then $\bigvee A\in\Ll$
 and $\bigwedge A\in\Ll$.
\end{enumerate}

Working for a moment in a possibly larger universe,
given $\varphi\in\Ll$ and a 
set $X\sub V_\beta^{M_\om}$, we define satisfaction 
$X\sats\varphi$ recursively in the obvious manner:
$X\sats P_x$ iff $x\in X$,
 $X\sats\neg\varphi$ iff $\neg (X\sats\varphi)$,
$X\sats\bigvee A$ iff
$X\sats\varphi$ for some $\varphi\in A$,
and $X\sats\bigwedge A$ iff $X\sats\varphi$ for all $\varphi\in A$.

Recall 
$x^*=i_{\om,\om+\om}(x)=\widetilde{i}_{0\om}^{M_\om}(x)$ in $M_0$.
Work again in $M_\om$. 
For $\varphi\in\Ll$, let
\[ \mathscr{E}_\varphi=\{X\sub V_{\beta^*}\bigm|X\sats\varphi^*\}.\]
Let $\Ll'$ be the class of formulas $\varphi\in\Ll$
such that $\mathscr{E}_\varphi\neq\emptyset$.
Given $\varphi,\psi\in\Ll'$, set $\varphi\equiv\psi$
iff $\mathscr{E}_\varphi=\mathscr{E}_\psi$.
So there are only set-many (in fact at most $V_{\beta^*+2}$-many)
equivalence classes. Let $\PP$ be the partial order
whose conditions are the equivalence classes $[\varphi]$,
and with
$[\varphi]\leq[\psi]\iff\mathscr{E}_{\varphi}\sub\mathscr{E}_{\psi}$.
This is a partial order (forcing equivalent to a partial order $\sub 
V_{\beta^*+2}$).

Now work in $V$ and let $X\sub V_\beta$. Let 
\[ G_X=\{[\varphi]\in\PP\bigm|\varphi\in\Ll'\text{ and }X\sats\varphi\}.\]
Here if $\varphi,\psi\in\Ll'$
and $\varphi\equiv\psi$ then $X\sats\varphi\Leftrightarrow X\sats\psi$.
For suppose that $X\sats\varphi\wedge\neg\psi$.
Let $n<\om$ be large enough that 
$i_{nm}(\varphi\wedge\neg\psi)=\varphi\wedge\neg\psi$
and $i_{nm}(\beta)=\beta$
for all $m\in[n,\om)$. Then note that
\[ V\sats\text{``There is }X'\sub V_\beta\text{ such that 
}X'\sats\varphi\wedge\neg\psi\text{''}\]
(as witnessed by $X$). Therefore, applying $i_{n\om}$, we have
\[ M_\om\sats\text{``There is }X'\sub V_{\beta^*}\text{ such that }
 X'\sats(\varphi\wedge\neg\psi)^*\text{''}.\]
But then for any such $X'$, we have $X'\sats\varphi^*$ and 
$\neg(X'\sats\psi^*)$ (as 
$(\varphi\wedge\neg\psi)^*=\varphi^*\wedge\neg(\psi^*)$), so 
$\mathscr{E}_\varphi\neq\mathscr{E}_\psi$,
a contradiction.

Now we claim that $G_X$ is $(M_\om,\PP)$-generic;
and then it is easy to see that $M_\om[G_X]=M_\om[X]$.
For $G_X$ is easily a filter, so we just need to verify genericity.
So let $D\in\pow(\PP)\inter M_\om$ be dense.
We must see that $G_X\inter D\neq\emptyset$.
Note that since $D$ is a set, for each $\xi\in\OR$, $M_\om$ has the set
\[ E=\{\varphi\in\Ll'\inter V_\xi\bigm|[\varphi]\in D\},\]
and taking $\xi$ large enough,
$D=\{[\varphi]\bigm|\varphi\in E\}$.
Now $\bigvee E\in\Ll$. We must see that $X\sats\bigvee E$,
because then $X\sats\varphi$ for some $\varphi\in E$,
and therefore $G_X\inter D\neq\emptyset$.

So suppose $X\sats\neg\bigvee E$.
Let $n<\om$ with $i_{n,n+1}(E,\beta)=(E,\beta)$.
Then
\[ V\sats\text{``There is }X'\sub V_\beta\text{ such that }X'\sats\neg\bigvee 
E\text{''},\]
\[ M_\om\sats\text{``There is }X'\sub V_{\beta^*}\text{ such that 
}X'\sats(\neg\bigvee E)^*\text{''}.\]
So $\neg\bigvee E\in\Ll'$. By the density of $D$,
there is $\varphi\in E$ such that $[\varphi]$ is compatible
with $[\neg\bigvee E]$. So let $\psi\in\Ll'$ with
$[\psi]\leq[\varphi]\text{ and }[\psi]\leq[\neg\bigvee E]$.
Then
$\emptyset\neq\mathscr{E}_\psi\sub\mathscr{E}_\varphi\inter\mathscr{E}_{
\neg\bigvee E}$.
So let $Y\in M_\om$, with $Y\sub V_{\beta^*}$, witness that 
$\mathscr{E}_\varphi\inter\mathscr{E}_{\neg\bigvee E}\neq\emptyset$. Then
$Y\sats\varphi^*\text{ and }
 Y\sats(\neg\bigvee E)^*$,
but note 
$(\neg\bigvee 
E)^*=\neg\bigvee(E^*)\equiv\bigwedge\{\neg\varrho\bigm|\varrho\in E^*\}$
and $\varphi\in E$, so $\varphi^*\in E^*$, so
$Y\sats\varphi^*\text{ and }Y\sats\neg\varphi^*$,
contradiction.
\end{proof}

We can immediately deduce:
\begin{cor}
 Suppose $(V,j)\sats\ZFR$. Let $\alpha<\om^2$.
Then every set in $V$ is set-generic over $M_\alpha$.
\end{cor}
\begin{proof}
 We have this for $\alpha<\om+\om$, by Theorem \ref{tm:every_om}.
 Consider then $\alpha=\om+\om$. By Theorem \ref{tm:every_om},
 applied \emph{in} $M_\om$, we get that every set in $M_\om$
 is set-generic over $M_{\om+\om}$. Therefore
 given any $\alpha\in\OR$, we can find a set-generic
 extension $M_{\om+\om}[G]$
 of $M_{\om+\om}$ such that $V_\alpha^{M_{\om+\om}[G]}=V_\alpha^{M_\om}$.
 But then for any set $X\in V$, by taking $\alpha$ high enough,
 we get that $X$ is set-generic over $M_{\om+\om}[G]$.
 So we have added $X$ by a 2-step iteration over $M_{\om+\om}$, which suffices. 
 Clearly this generalizes to all $\alpha<\om+\om$.
\end{proof}

We next generalize these initial facts to arbitrary iterates.
\begin{dfn}
For a limit $\alpha$, say  $x$ is \emph{$\alpha$-eventually stable} 
(\emph{$\alpha$-es})
 iff
$x\in\bigcap_{\xi<\alpha}M_\xi$
and there is $\xi<\alpha$ such that
$i_{\xi\gamma}(x)=x\text{ for all }\gamma\in[\xi,\alpha)$.
Say that $x$ is \emph{hereditarily $\alpha$-eventually stable}
(\emph{$\alpha$-hes})
iff every $y\in\trancl(\{x\})$ is $\alpha$-es.
 For $\alpha$-es $x$, let
 $x^*_\alpha=\lim_{\xi<\alpha}i_{\xi\alpha}(x)$.
 Note here that $i_{\xi\alpha}(x)=i_{\gamma\alpha}(x)$
 whenever $\xi\leq\gamma<\alpha$ and $\xi$ witnesses the $\alpha$-eventual 
stability of $x$.\end{dfn}
\begin{dfn}
Given an ordinal $\alpha>0$, the \emph{eventual ordertype}
$\eot(\alpha)$ of $\alpha$ is the least $\chi$
such that $\alpha=\xi+\chi$ for some $\xi<\alpha$.
\end{dfn}

Note that (i) $\alpha$ is a successor iff $\eot(\alpha)=1$ iff $\eot(\alpha)$ 
is a successor; (ii) $\eot(\alpha)=\eot(\eot(\alpha))$;
(iii) $\alpha=\xi+\eot(\alpha)$
for all sufficiently large $\xi<\alpha$.
Clearly $\eta=\eot(\alpha)\leq\alpha\leq\kappa_\alpha$,
and if $\eta<\kappa_\alpha$ 
then $\eta^*_\alpha=\eta$,
whereas if $\eta=\alpha=\kappa_\alpha$ then 
$\alpha<\alpha^*_\alpha=\kappa_{\alpha^*_\alpha}$:

\begin{lem}\label{lem:star_general}Suppose $(V,j)\sats\ZFR$. Let $\alpha$ be a 
limit ordinal.
 Then:
 \begin{enumerate}
  \item $N_\alpha=\{x\bigm|x\text{ is }\alpha\text{-hes}\}$
\item 
$x^*_\alpha=i_{\alpha,\alpha+\eta^*_\alpha}(x)=
\widetilde{i}^{M_\alpha}_{0\eta^*_\alpha} (x)$
where $\eta=\eot(\alpha)$.
 \end{enumerate}
\end{lem}
\begin{proof}
First note that all ordinals are $\alpha$-es, since $M_\alpha$ is wellfounded.

Now let $x\in M_\alpha$; we show that $x$ is $\alpha$-hes.
We have $M_\alpha\sub M_\xi$ for each $\xi<\alpha$.
Assume by induction that every $y\in\trancl(x)$ is $\alpha$-es.
So we  need to see  $x$ is $\alpha$-es.
Let $\xi<\alpha$ be such that $x\in\rg(i_{\xi\alpha})$
and $\xi+\eta=\alpha$ where $\eta=\eot(\alpha)$,
and $i_{\xi\gamma}(\eta)=\eta$ for all $\gamma\in[\xi,\alpha)$.
We claim that $i_{\xi\gamma}(x)=x$ for all $\gamma\in[\xi,\alpha)$.
For
$x=i_{\xi\alpha}(\bar{x})=\widetilde{i}^{M_\xi}_{0\eta}(\bar{x})$
for some $\bar{x}$,
so 
$i_{\xi\gamma}(x)=i_{\xi\gamma}(\widetilde{i}^{M_\xi}_{0\eta})(i_{\xi\gamma}
(\bar{x}))=\widetilde{i}^{M_\gamma}_{0\eta}(i_{\xi\gamma}(\bar{x}
))$,
but $\gamma+\eta=\alpha$, so 
$\widetilde{i}^{M_\gamma}_{0\eta}(i_{\xi\gamma}(\bar{x}))=
i_{\gamma\alpha}(i_{\xi\gamma}(\bar{x}))=
i_{\xi\alpha}(\bar{x})=x$,
as desired.

Now let $x\in M_\alpha$ and $\eta=\eot(\alpha)$.
We must see that $x^*_\alpha=i_{\alpha,\alpha+\eta^*_\alpha}(x)$.
Let $\xi<\alpha$ be such that
(i) $\xi+\eta=\alpha$,
(ii) $x\in\rg(i_{\xi\alpha})$,
and (iii) $i_{\xi\gamma}(x,\eta)=(x,\eta)$ for all $\gamma\in[\xi,\alpha)$.
Then $(x^*_\alpha,\eta^*_\alpha)=i_{\xi\alpha}(x,\eta)$;
let $\bar{x}$ be such that $i_{\xi\alpha}(\bar{x})=x$.
Since $M_\xi\sats$``$x=\widetilde{i}_{0\eta}(\bar{x})$'',
applying $i_{\xi\alpha}$, we get
\[ M_\alpha\sats\text{``}x^*_\alpha=\widetilde{i}_{0\eta^*_\alpha}(x)\text{''}, 
\]
so $i_{\alpha,\alpha+\eta^*_\alpha}(x)=x^*_\alpha$, as desired.
\end{proof}

\begin{tm}\label{tm:every_set_generic_over_N_tau}
 Suppose $(V,j)\sats\ZFR$ and let $\tau\in\OR$.
 Then every set in $V$ is set-generic over $M_\tau$.
\end{tm}
\begin{proof}
 By induction on $\tau$. Successor steps are trivial,
 so let $\tau$ be a limit and $\chi=\eot(\tau)$. We write 
$x^*=x^*_\tau=i_{\tau,\chi^*_\tau}(x)$ throughout (as computed in $(V,j)$).

Let $\alpha\in\OR$
 and $\beta=\sup i_{0\tau}``\alpha$. Let $X=i_{0\tau}``V_\alpha$.
 It suffices to see that $X$ is set-generic over $M_\tau$.
 We will just use the fact that $X\sub V_{\beta}^{M_\tau}$.

Now $X$ is set-generic over $M_\xi$ for each $\xi<\tau$.
Fix $\delta\in\OR$
such that for each $\xi<\tau$, there is an 
$(M_\xi,\Coll(\om,V_\delta^{M_\xi}))$-generic
 $G$ such that $X\in M_\xi[G]$.

 Work in $M_\tau$. Let $\Ll$ be as before
 (with $P_x$ for each $x\in V_\beta$). 
 Let $\QQ=\Coll(\om,V_{\delta^*})$.
For $\varphi\in\Ll$,  let $\mathscr{E}_{\varphi}$ be the $\QQ$-name
  for the set of all 
$X\in\pow(V_{\beta^*}^{\widetilde{M}_{\chi^*}})\inter V^\QQ$
such that $X\sats\varphi^*$.
 Let $\Ll'$ be the class of all $\varphi\in\Ll$
 such that $\QQ\forces\mathscr{E}_\varphi\neq\emptyset$
 (equivalently, there is $p\in\QQ$ such that 
$p\forces\mathscr{E}_\varphi\neq\emptyset$).
Given $\varphi,\psi\in\Ll'$, let  $\varphi\equiv\psi$ iff 
$\QQ\forces\mathscr{E}_\varphi=\mathscr{E}_\psi$.
 Let $[\varphi]$ be the equivalence class of $\varphi\in\Ll'$.
 Let $\PP$ be the partial order whose conditions
 are the equivalence classes $[\varphi]$,
 with 
$[\varphi]\leq[\psi]$  iff $\QQ\forces\mathscr{E}_\varphi\sub\mathscr{E}_\psi$.

Now work in $V$. Let $G_X=\{\varphi\in\Ll'\bigm|X\sats\varphi\}$.
We claim that $G$ is $(M_\tau,\PP)$-generic;
it follows that $M_\tau[G_X]=M_\tau[X]$.
This is proven like before. Let $D\in\pow(\PP)\inter M_\tau$ be dense.
Let $E\in M_\tau$ be a such that $D=\{[\varphi]\bigm|\varphi\in E\}$.
We need to see that $X\sats\bigvee E$, so suppose otherwise.
Then for all sufficiently large $\xi<\tau$,
\[ M_\xi\sats\text{let }Y=V_\beta^{\widetilde{M}_\chi};\text{ then 
}\Coll(\om,V_\delta)\forces\exists X'\sub\check{Y}\text{ such that 
}X'\sats\neg\bigvee E,\]
(as witnessed by $X$, by induction and choice of $\delta$), which 
gives that
\[ M_\tau\sats\text{``let 
}Y=V_{\beta^*}^{\widetilde{M}_{\chi^*}};\text{ then }\QQ\forces\exists 
X'\sub\check{Y}\text{ such that }X'\sats(\neg\bigvee E)^*.\]
So $(\neg\bigvee E)\in\Ll'$, so $[\neg\bigvee E]\in\PP$.

Work in $M_\tau$. By the density of $D$ and since $[\neg\bigvee E]\in\PP$,
there is $\varphi\in E$ such that $[\varphi]$ is compatible
with $[\neg\bigvee E]$. Let $[\psi]\in\PP$
with $[\psi]\leq[\varphi]$ and $[\psi]\leq[\neg\bigvee E]$.
Let $G$ be $\Coll(\om,V_{\delta^*})$-generic.
Then in $V[G]$, there is $X'$ such that $X'\sats\psi^*$,
and by the definition of $\leq$, therefore $X'\sats\varphi^*$
and $X'\sats(\neg\bigvee E)^*$, so $X'\sats\bigwedge_{\varrho\in 
E^*}\neg\varrho$, so $X'\sats\neg\varphi^*$,
a contradiction.
\end{proof}

We deduce the existence of some interesting forcing extensions under 
$\ZFR$. Recall the notion of \emph{virtual large cardinals} from 
\cite{virtual_lc}.
\begin{dfn}\label{dfn:virtually_Berkeley}
$\delta$ is \emph{virtually Berkeley}
 iff for every transitive $M$ with $\delta\in M$
 and every $\eta<\delta$, there is a forcing $\PP$ which forces
 the existence of an elementary $j:M\to M$ with $\eta<\crit(j)<\delta$.
 Say that $\delta$ is \emph{$V_\alpha$-preserving-virtually Berkeley}
 iff $\delta$ is virtually Berkeley, as witnessed by forcings
 $\PP$ which preserve $V_\alpha$.
\end{dfn}

\begin{cor}\label{cor:forcing_extensions}
 Suppose $(V,j)\sats\ZFR$. Let $\kappa=\crit(j)$ and $\lambda=\kappa_{\om}(j)$. 
Then there are partial orders $\PP_0,\PP_1\in V$ such that:
\begin{enumerate}
 \item\label{item:virtually_Berkeley} $\kappa$ is 
$V_\kappa$-preserving-virtually Berkeley.
  \item\label{item:PP_0} $\PP_0\forces V_{\check{\kappa}}=\check{V_\kappa}$
  and there is an elementary 
 $k:V_{\check{\kappa}}\to V_{\check{\kappa}}$
  such that
  \begin{enumerate}[label=\tu{(}\alph*\tu{)}]
\item $\kappa_{\om}(k)=\check{\kappa}$ \tu{(}hence 
$\cof(\check{\kappa})=\om$\tu{)}
and
\item $(\check{V_\lambda},\check{j}\rest\check{V_\lambda})$ is 
the $\om$th 
iterate of $(\check{V_\kappa},k)$.
\end{enumerate}
\item\label{item:PP_1} for each $\beta<\kappa$ there is 
$\bar{\lambda}\in(\beta,\kappa)$
such that $\PP_1$ forces that 
$V_{\check{\bar{\lambda}}}=\check{V_{\overline{\lambda}}}$
and $(\bar{\lambda}^+)=\check{\kappa}$
and there is a rank-into-rank embedding $k:V_{\check{\bar{\lambda}}}\to 
V_{\check{\bar{\lambda}}}$ such that $(\check{V_\lambda},\check{j}\rest 
\check{V_\lambda})$
is the $\check{\kappa}$th iterate of $(V_{\check{\bar{\lambda}}},k)$.
 \end{enumerate} 
\end{cor}
\begin{proof}
 The corresponding things are forceable over $M_\om$
 (for parts \ref{item:virtually_Berkeley} and \ref{item:PP_0})
 and over $M_{(\lambda^+)}$ (for part \ref{item:PP_1}) with respect to 
$\lambda$,
 respectively. So the elementarity of $i_{0\om}$
 and $i_{0(\lambda^+)}$ pulls  them  back to $M_0=(V,j)$.
\end{proof}

\subsection{Intersections of iterates}

We consider now intersections of the form $\bigcap_{\xi<\tau}N_\xi$.

\begin{dfn}
Suppose $(V,j)\sats\ZFR$ and let $\tau$ be a limit ordinal
and $\beta<\tau$. Let $x\in\bigcap_{\xi<\tau}N_\xi$.
Then $x$ is \emph{$(\beta,\tau)$-stable}
iff $i_{\beta\gamma}(x)=x$ for all $\gamma\in[\beta,\tau)$.
\end{dfn}

\begin{tm}
Suppose $(V,j)\sats\ZFR$ and let $\tau$ be a limit ordinal
 and $\chi=\eot(\tau)$.
 Let $J=\bigcap_{\xi<\tau}N_\xi$.
 Let $\xi<\tau$ be such that $\chi$ is $(\xi,\tau)$-stable
 and either
 \begin{enumerate}[label=\tu{(}\roman*\tu{)}]
  \item\label{case:eta<crit}
$\chi<\crit(j_\xi)$,
 or
 \item\label{case:eta=sup}
$\kappa_\xi<\chi=\tau=\kappa_\tau=\sup_{\gamma<\tau}\kappa_\gamma$.
\end{enumerate}
 Let $\mu=\cof^{M_\xi}(\tau)=\cof^{M_\xi}(\chi)$, so 
$i_{\xi\gamma}(\mu)=\cof^{M_\gamma}(\chi)$
 for all $\gamma\in[\xi,\tau)$. Then:
 \begin{enumerate}
  \item\label{item:J_sats_ZF} $J$ is transitive, proper class
  with $N_\tau\sub J$ and $(J,M_\tau)\sats\ZF$.
  \item\label{item:mu>om} if $\mu>\om$ then $N_\tau=J$.
  \item\label{item:mu=om} if $\mu=\om$ then $N_\tau\psub 
J\sats$``$\cof(\gamma)=\om$
  for all
$N_\tau$-regular 
$\gamma\in[\kappa_\tau,\kappa_{\tau+\om})$''.
 \end{enumerate}
\end{tm}
\begin{proof}
 Part \ref{item:J_sats_ZF}: $J$ is clearly  transitive with  
$N_\tau\sub J$.
 And $(J,M_\tau)$ is a class of $M_\xi=(N_\xi,j_\xi)$ for each 
$\xi<\tau$.
 So $J$ satisfies Powerset: if $X\in J$ and $\xi<\tau$ then because $J$ is a 
class of $M_\xi$,
 we have $\pow(X)\inter J\in M_\xi$; therefore $\pow(X)\inter J\in J$.
Likewise, for each $\alpha\in\OR$, we have $V_\alpha\inter J\in J$
and $\left<V_\beta\inter J\right>_{\beta<\alpha}\in J$.
Separation (with respect to the predicate $M_\tau$) is similar to Powerset, 
and Collection follows
from the preceding remarks.
 
 Part \ref{item:mu>om}: Suppose $\mu>\om$.
 Let $x\in J$. By Lemma \ref{lem:star_general} and transitivity, it suffices 
to see that $x$ is 
$\tau$-es.  Suppose not.  Then for each $\gamma<\tau$ there is 
$\gamma'\in(\gamma,\tau)$
such that $i_{\gamma\gamma'}(x)\neq x$.

Work in $M_\xi$, where $\cof(\chi)=\mu>\om$.
We have $\xi+\chi=\tau$.
 Define a sequence $\left<\xi_n\right>_{n<\om}$ as follows:
 Set $\xi_0=0$. Given $\xi_n$, let $\xi_{n+1}$
 be the least $\xi'\in(\xi_n,\eta)$ such that 
$\widetilde{i}^{M_\xi}_{\xi_n\xi'}(x)\neq x$.
Then since $\cof(\chi)>\om$, we have $(\sup_{n<\om}\xi_n)=\xi_\om<\chi$.
So $x\in\widetilde{N}_{\xi_\om}$, so by Lemma \ref{lem:star_general},
$x$ is $\xi_\om$-es, a contradiction.

 Part \ref{item:mu=om}: Suppose $\mu=\om$.
 Work in $M_\xi$. Let 
$\left<\xi_n\right>_{n<\om}$ 
with $\sup_{ n<\om}\xi_n=\chi$ and $\xi_0=0$. Define
$\left<\xi'_n\right>_{1\leq n<\om}$ by 
$\widetilde{i}^{M_\xi}_{\xi_n}(\xi_n)$ and $\xi'_0=0$.
Note 
$\left<\xi'_n\right>_{n>m}\in M_{\xi'_m}$
(by commutativity of the maps), so in fact
$\left<\xi'_n\right>_{n<\om}\in J$.
Note $\xi_n\leq\xi'_n<\xi'_{n+1}<\chi$ (since
 $\chi$ is $(\xi,\tau)$-stable)
and $\sup_{n<\om}\xi'_n=\sup_{n<\om}\xi_n=\chi$.

Renaming, we may assume $\left<\xi_n\right>_{n<\om}\in J$.
Now let $\gamma\in[\kappa_\tau,\kappa_{\tau+\om})$
be regular in $N_\tau$. We need
$J\sats$``$\cof(\gamma)=\om$''.
We may assume $\gamma\in\rg(i_{\xi\tau})$;
let $i_{\xi\tau}(\bar{\gamma})=\gamma$.
So $M_\xi\sats$``$\bar{\gamma}$ is regular''.
Let
$\bar{\gamma}_n=i_{\xi,\xi+\xi_n}(\bar{\gamma})$,
so $i_{\xi+\xi_n,\tau}(\bar{\gamma}_n)=\gamma$.
Let
\[ \bar{\zeta}_{n+1}=\sup i_{\xi+\xi_n,\xi+\xi_{n+1}}``\bar{\gamma}_n.\]
Now $\bar{\zeta}_{n+1}<\bar{\gamma}_{n+1}$. 
For
$j_{\xi+\xi_n}(\gamma')>\gamma'$,
for all 
$\gamma'\in[\kappa_{\xi+\xi_n},\kappa_{\xi+\xi_n+\om})$.
But $j_{\xi+\xi_n}:M_{\xi+\xi_n}\to M_{\xi+\xi_n}$,
so by the regularity of $j_{\xi+\xi_n}(\bar{\gamma}_n)$ in $M_{\xi+\xi_n}$,
and since $j_{\xi+\xi_n}\rest\bar{\gamma}_n\in M_{\xi+\xi_n}$, we get
\[ (\sup j_{\xi+\xi_n}``\bar{\gamma}_n)<j_{\xi+\xi_n}(\bar{\gamma}_n),\]
and by commutativity, therefore
$\bar{\zeta}_{n+1}<\bar{\gamma}_{n+1}$.
Let $\zeta_{n+1}=i_{\xi+\xi_{n+1},\tau}(\bar{\zeta}_{n+1})$.
Then $\sup_{n<\om}\zeta_{n+1}=\gamma$,
and  note $\left<\zeta_{n+1}\right>_{n<\om}\in J$,
using $\left<\xi_n\right>_{n<\om}\in J$.
\end{proof}

\subsection{The (infinite) iterates $N_\tau$ are not set grounds of $V$}

In contrast to the fact that every set is generic over $N_\tau$, we have:

\begin{tm}\label{tm:N_tau_not_set-ground}
 Suppose $(V,j)\sats\ZFR$.  Then:
 \begin{enumerate}
 \item\label{item:not_proper_ground} For all $\tau\in\OR$ with 
$\tau\geq\om$, 
$N_\tau$ is not a  set ground of 
$V$.
  \item\label{item:j_rest} For $\xi<\tau\in\OR$, we have:
  \begin{enumerate}
\item\label{item:rest_in}  $i_{\xi\tau}\rest V_\alpha^{N_\xi}\in N_\tau$
  for all $\alpha<\kappa_{\xi+\om}$.
  \item\label{item:rest_not_in}
If $\xi+\om\leq\tau$ then $i_{\xi\tau}\rest\kappa_{\xi+\om}\notin 
N_\tau$.
\end{enumerate}
 \end{enumerate}
\end{tm}
\begin{proof}
Part \ref{item:rest_in}: If $\tau<\xi+\om$ this is immediate 
because
$N_\tau=N_\xi$. So suppose $\xi+\om\leq\tau$.
Let $n<\om$ be such that $\alpha<\kappa_{\xi+n}$.
Then $k=i_{\xi,\xi+n}\rest V_\alpha^{N_\xi}\in N_{\xi+n}=N_\xi$,
and note $i_{\xi\tau}\rest 
V_\alpha^{N_\xi}=i_{\xi+n,\tau}\com k=i_{\xi+n,\tau}(k)\in N_\tau$.

Part \ref{item:rest_not_in}:
Suppose first $\tau=\xi+\om$.
Then $i_{\xi,\xi+\om}``\kappa_{\xi+\om}$ is cofinal in $\kappa_{\xi+\om+\om}$,
and $\cof^{N_{\xi+\om}}(\kappa_{\xi+\om})=\kappa_{\xi+\om}=\crit(j_{\xi+\om})$.
But $\cof^{N_{\xi+\om}}(\kappa_{\xi+\om+\om})=\om$, because 
$\kappa_{\xi+\om+\om}=\widetilde{\kappa}_\om^{M_{\xi+\om}}$.
Therefore $i_{\xi,\xi+\om}\rest\kappa_{\xi+\om}\notin N_{\xi+\om}$.

Now suppose $\tau>\xi+\om$. Since $i_{\xi+\om,\tau}$ is computed by 
$M_{\xi+\om}$, and by commutativity, if $i_{\xi\tau}\rest\kappa_{\xi+\om}\in 
N_\tau\sub N_{\xi+\om}$, then $i_{\xi,\xi+\om}\rest\kappa_{\xi+\om}\in 
N_{\xi+\om}$, a contradiction.
 
Part \ref{item:not_proper_ground}: 
By part \ref{item:j_rest},
$i_{0\tau}\rest\kappa_\om$ is a set of   ordinals in $V\cut N_\tau$,
and $i_{0\tau}:V\to N_\tau$ is elementary.
So by Theorem \ref{tm:no_gen_emb_into_V_with_new_O_set},  $N_\tau$ is not a 
set-ground of $V$.
\end{proof}

\subsection{Non-definability of $j\rest\OR$}

\begin{rem}\label{rem:j_rest_OR_not_def_from_params}
Recall Suzuki \cite{suzuki_no_def_j} showed that no elementary $j:V\to V$ 
is definable
from parameters. Variants  of this are proved in 
\cite{gen_kunen_incon}, \cite{cumulative_periodicity}, \cite{ZF_extenders}.
In particular, \cite[Theorem 35]{gen_kunen_incon}
establishes the stronger fact that no elementary $k:\HOD\to\HOD$ is definable 
from parameters (hence if $j:V\to V$ is elementary then $j\rest\OR$ is not 
definable from parameters). We sketch now an alternate proof of this fact, 
closer to the methods of this paper.
Suppose $k:\HOD\to\HOD$ is definable from parameter $x$. Using Vopenka (or 
maybe a slight adaptation thereof 
to $\ZF$, if it is not standard; see Fact \ref{fact:Vopenka}),
there is a $\HOD$-set-generic filter $G_x$ such that
$\HOD[G_x]=\HOD_x$.\footnote{We cannot expect that $x\in\HOD[G_x]$ in general,
since $\HOD[G_x]\sats\ZFC$, and $x$ is arbitrary. However,
the usual arguments still give $\HOD[G_x]=\HOD_x$.}
So $k\rest\HOD$ is amenable to $\HOD[G_x]$.
But this is impossible by \cite[Corollary 9]{gen_kunen_incon};
or alternatively, using Lemma  \ref{lem:no_gen_emb_into_V_with_new_O_set}, we 
get that $k\rest\HOD$
is amenable to $\HOD$, contradicting Kunen.

In this section we strengthen the fact that $j\rest\OR$ is not definable
from parameters, if $j:V\to V$ is elementary, 
showing $j\rest\OR$ cannot be definable 
from 
parameters over the structure $(V,M_\om)$,
that is, allowing $M_\om$ as a predicate.

Also, it is shown in \cite[Theorem 12]{gen_kunen_incon} under 
$\ZFC_2$ that if  $j:V\to M$ is an elementary 
 with $M$ wellfounded, then $V=\Hull^V(\rg(j))$.
 But the proof relies strongly on $\AC$.
 It trivially fails with respect to $j$
 when $(V,j)\sats\ZFR$ (we have $\rg(j)=\Hull^V(\rg(j))$ by elementarity).
 It also fails with respect to $j=i_{0\om}$,
  as shown in the following theorem (but since $N_\om\neq V$,
 this time it is not trivial).
 Note that $\Hull^M(X)$ was defined in \S\ref{subsec:notation}.
\end{rem}

\begin{tm}\label{tm:V=HOD(X)_no_R}
Assume $(V,j)\sats\ZFR$. Then:
\begin{enumerate}
 \item\label{item:V_not_Hull} For no set $X$ is $V=\Hull^{(V,M_\om)}(N_\om\cup 
X)$.
\item\label{item:j_rest_OR_not_def} $j\rest\OR$ and $i_{0\om}\rest\OR$ are not 
definable from parameters over $(V,M_\om)$.
\end{enumerate}
\end{tm}

\begin{proof}
Part \ref{item:V_not_Hull}:
Suppose $X\in V$ and
 $V=\Hull^{(V,M_\om)}(N_\om\cup X)$. Let 
$k_n=i_{n\om}\rest X^{<\om}$, for $n<\om$.
 So $\rg(k_n)\sub N_\om$.
 Let $\vec{k}=\left<k_n\right>_{n<\om}$.
 Then $(X,\vec{k})\in V$, so by Theorem \ref{tm:every_om}, we have 
$(X,\vec{k})\in 
N_\om[G]$ for some set-generic extension $N_\om[G]$ of $N_\om$.
We will show that $N_\om[G]=V$, contradicting Theorem 
\ref{tm:N_tau_not_set-ground}.

Let $Y\in V$. By our contradictory assumption,
there is $\alpha\in\OR$ and 
$n<\om$ such that
$Y$ is encoded
into $t=\Th_{\Sigma_n}^{(V,M_\om)}((N_\om\inter V_\alpha)\cup X)$. So it 
suffices to see
$t\in N_\om[G]$. We have $(N_\om\inter V_\alpha)\cup X\in N_\om[G]$,
so we just need to verify the appropriate instance of Separation.

So let $\varphi$ be $\Sigma_n$ and $\vec{y}\in(N_\om\inter 
V_\alpha)^{<\om}$
and $\vec{x}\in X^{<\om}$. For each
$m<\om$, note
\[ 
[(V,M_\om)=(N_m,M_\om)\sats\varphi(\vec{y},\vec{x})]\iff 
[(N_\om,M_{\om+\om})\sats\varphi(i_{m\om}(\vec{y}),i_{m\om}(\vec{x}))].\]
But $i_{m\om}(\vec{x})=k_m(\vec{x})$ for all $m<\om$,
and $i_{m\om}(\vec{y})=\vec{y}^*=i_{\om,\om+\om}(\vec{y})$ for all 
sufficiently large $m<\om$.
Therefore
\[ [(V,M_\om)\sats\varphi(\vec{y},\vec{x})]\iff
\exists \ell<\om\all m\geq\ell\ [(N_\om,M_{\om+\om})\sats 
\varphi(\vec{y}^*,k_m(\vec{x}))].\]
Since $\left<k_m\right>_{m<\om}\in N_\om[G]$ 
and $i_{\om,\om+\om}\rest(N_\om\inter V_\alpha)\in N_\om$,
we get $t\in N_\om[G]$.\footnote{\label{ftn:ground_def_claim}Of 
course,
in $N_\om[G]$ we can compute the relevant truth of $(N_\om,M_{\om+\om})$,
just noting that whenever
 $(V_\xi,V_\xi\inter M_\om)\elem_n (V,M_\om)$, then 
$(V^{N_\om}_{i_{0\om}(\xi)},V^{N_\om}_{i_{0\om}(\xi)}\inter M_{\om+\om})\elem_n 
(N_\om,M_{\om+\om})$. (In the notes \cite{reinhardt_non-definability},
at this point it was claimed that we could alternatively
use the definability of set-grounds here to allow $N_\om[G]$ to compute $N_\om$,
but this was an oversight, since $\AC$ fails, and it is not clear
to the author whether one can get around that failure here.)}

Part \ref{item:j_rest_OR_not_def}: Arguing as above, note that for each $x\in 
V$, letting $x_n=i_{n\om}(x)$, and letting $N_\om[G]$ be a set-generic
extension such that $\left<x_n\right>_{n<\om}\in N_\om[G]$,
we get that every set $A\sub N_\om$ which is definable over $(V,M_\om)$
from parameters in $N_\om\cup\{x\}$, is in $N_\om[G]$.

Now suppose
$i_{0\om}\rest\OR$ is definable over $(V,M_\om)$ from elements of 
$N_\om\cup\{x\}$.
Then by the previous paragraph, $i_{0\om}\rest\OR$ is amenable to $M_\om[G]$.  
But then much as
 in the proof of Lemma \ref{lem:no_gen_emb_into_V_with_new_O_set}, we can show 
that $i_{0\om}\rest\OR$
 is amenable to $M_\om$, contradicting Theorem \ref{tm:N_tau_not_set-ground}.
 That is, fix $\alpha\in\OR$.
We can find an $X\in M_\om$ and $\beta\in\OR$ with $X\sub\beta$ and $X$ of 
ordertype $\alpha$, such that $i_{0\om}\rest\alpha\in M_\om$. Let 
$\pi:\alpha\to X$ be the increasing enumeration. Then 
$i_{0\om}(\pi),i_{0\om}(X)\in N_{\om+\om}\sub N_\om$, and 
$i_{0\om}(\pi):i_{0\om}(\alpha)\to i_{0\om}(X)$
is the increasing enumeration of $i_{0\om}(X)$.
Since $i_{0\om}\com\pi=i_{0\om}(\pi)\com i_{0\om}\rest\alpha$,
and $i_{0\om}\rest X=\rg(\pi)\in N_\om$, we therefore get 
$i_{0\om}\rest\alpha\in N_\om$, as desired.

Finally suppose  $j\rest\OR$ is definable over $(V,M_\om)$ from elements of 
$N_\om\cup\{x\}$.\footnote{Note this doesn't immediately contradict
the previous paragraph, because from $j\rest\OR$ alone it doesn't seem that
one can directly compute $i_{0\om}\rest\OR$. We can compute $i_{0n}\rest\OR$
for $n<\om$, since this is the same as $j\com\ldots\com j\rest\OR$. But
for $i_{0\om}\rest\OR$ it seems we need the direct limit of all $j_n\rest\OR$.}
Note $j\rest\kappa_\om\notin N_\om$, since
otherwise $\cof^{N_\om}(\kappa_\om)=\om$. But $j\rest N_\om:M_\om\to M_\om$  
is elementary. So we can just use essentially the same argument
as in the previous paragraph for a contradiction (but slightly
easier).
\end{proof}

\begin{rem}\label{rem:goldberg_usuba}
An immediate corollary is that there is no $X\in V$
 such that $V=\HOD(X)$.
And if $(V_\delta,V_{\delta+1})\sats\ZF_2+$``$\delta$ is total 
Reinhardt'',
then there is no  $A\sub V_\delta$ and  $X\in V_\delta$ such that 
$V_\delta=(\HOD_A(X))^{(V_\delta,A)}$. For
 just let $j:(V_\delta,(A,X))\to (V_\delta,(A,X))$ be elementary,
 and run the obvious variant of the preceding proof.

But suppose $(V,j)\sats\ZFR$. In the notes
  \cite[v1 on arxiv.org]{reinhardt_non-definability}, the author asked
  whether possibly $V=\Hull^V(\OR\cup\rg(j))$,
 or $V=\Hull^V(M_\om\cup\rg(j))$,
 and whether $\AC$ can be set-forceable over $V$.
 These questions have since been answered:
 Goldberg \cite{goldberg_email} and Usuba \cite{usuba_email} have independently
 shown that $\AC$ is not 
set-forceable (not by considering the iterates $M_\alpha$,
 but by a more direct combinatorial method). In particular, there
 is no class $A$ and set $X$ such that $V=\HOD_A(X)$. Goldberg's proof
 also gives that if there is a Berkeley cardinal then $\AC$ is not 
set-forceable, and he also showed that
(i) $V=\Hull^V(M_\om\cup\rg(j))$, but
(ii) $V\neq\Hull^V(\OR\cup\rg(j))$. 

By an easier argument than that for (i),
we have $\Hull^{(V,j)}(N_\om)=V$, in contrast
to the fact that $\Hull^{(V,M_\om)}(N_\om)=N_\om$.
For $i_{0\om}$
is computable from $j$, and given $x\in V$,
we have $x=(i_{0\om})^{-1}(x')$ where $x'=i_{0\om}(x)\in N_\om$.
\end{rem}

\begin{rem}\label{rem:Vop}
Assume $\ZF$. Let $x\sub V_\alpha$.
We recall the definition of the 
\emph{Vopenka algebra}
$\Vop$ for 
 extending $\HOD$ generically to $\HOD_x=\HOD(\{x\})$.  Define $\Vop'$ in $V$ 
as follows.
 Conditions are non-empty $\OD$ subsets of $V_{\alpha+1}$;
 the ordering is $A\leq B\iff A\sub B$.
 Then $\Vop\in\HOD$ is the natural coding of $\Vop'$
 as a set of ordinals. Let $\pi:\Vop\to\Vop'$ be the natural isomorphism.
Now (with $x\sub V_\alpha$), define the filter
$G_x=\{p\in\Vop\bigm|x\in\pi(p)\}$.
\end{rem}
Then by the $\ZFC$ arguments, we have:
\begin{fact}\label{fact:Vopenka}
$G_x$ is $(\HOD,\Vop)$-generic and $\HOD[G_x]=\HOD_x$.\footnote{
The fact that $x\in\HOD[G_x]$, in the case that $x\sub\gamma\in\OR$,
is of 
course due to Vopenka. The author is not sure who first proved the 
entire result. Let $\varphi$ be a formula, $\eta,\gamma\in\OR$
and
$A=\{\xi<\gamma\bigm|\varphi(x,\xi,\eta)\}$.
We need $A\in\HOD[G_x]$.
For $\xi<\gamma$ let
$t'_\xi=\{y\bigm|y\in V_{\alpha+1}\text{ and 
}\varphi(y,\xi,\eta)\}$.
So $t'_\xi\in\Vop'$.
Let $\pi(t_\xi)=t'_\xi$.
Note $\left<t_\xi\right>_{\xi<\gamma}\in\HOD$,
and
$t_\xi\in G_x$ iff $x\in t'_\xi$ iff $\xi\in X$,
so $X\in\HOD[G_x]$.}
\end{fact}

Although every set is set-generic over $N_\tau$ for every $\tau<\OR$,
this fails at $\tau=\OR$, and although $\HOD\sub N_\om$,
there is $\tau<\OR$ with $\HOD\not\sub N_\tau$:
\begin{tm}
 Suppose $(V,j)\sats\ZFR$. Let $J=\bigcap_{\alpha\in\OR}N_\alpha$.
 Then:
 \begin{enumerate}
  \item\label{item:J_sats_ZF_again} $J$ is a transitive proper class and
$V^{N_\tau}_{\kappa_{\tau}}\elem V_{\kappa_{\tau+\om}}^{N_\tau}\elem J$ for 
each $\tau\in\OR$.
 \item\label{item:i_alpha,beta_J_to_J} For all $\alpha\leq\beta<\OR$,
 $i_{\alpha\beta}\rest J:J\to J$ is elementary.
 \item\label{item:j_rest_not_gen} 
 $j\rest V_{\kappa_\om}^{J}:V_{\kappa_\om}^{J}\to 
V_{\kappa_{\om}}^{J}$ is not set-generic over 
$J$.
\item\label{item:HOD_not_sub_J} $\HOD\not\sub J$.
 \end{enumerate}
\end{tm}
\begin{proof}
Part \ref{item:J_sats_ZF_again} is straightforward,
noting in particular that
$V_{\kappa_{\tau+\om}}^{N_\tau}=V_{\kappa_{\tau+\om}}^{N_{\tau+\om}}\sub J$.

Part \ref{item:i_alpha,beta_J_to_J}: Note that $J$ is defined in the same manner
over $M_\alpha=(N_\alpha,j_\alpha)$ as over $M_\beta$,
so $i_{\alpha\beta}(J)=J$.

Part \ref{item:j_rest_not_gen}: 
Let 
$R=(V_{\kappa_\om},j\rest V_{\kappa_\om})$.
Note  $R$ is iterable with  iterates
$R_\alpha=(V_{\kappa_{\alpha+\om}}^{M_\alpha},j_\alpha\rest 
V_{\kappa_{\alpha+\om}}^{M_\alpha})=(V_{\kappa_{\alpha+\om}}^J,j_\alpha\rest 
V_{\kappa_{\alpha+\om}}^J)$. But then if $R\in J[G]$,
we can form the iteration $\left<R_\alpha\right>_{\alpha\in\OR}$ there,
and considering some sufficiently large limit $\alpha$,
this collapses all cardinals in $(\kappa_\alpha,\kappa_{\alpha+\om}]$
in $J[G]$, although each $\kappa_{\alpha+n+1}$ is inaccessible
in $J$ and hence also in $J[G]$, a contradiction.

Part \ref{item:HOD_not_sub_J}: Let $G$ be 
$(V,\Coll(\om,V_{\kappa_\om}))$-generic.
Let $H=\HOD^{V[G]}$. By the homogeneity of the collapse, $H\sub\HOD$.
Let $x\in V[G]$ be a real coding $R$.
Then by Fact \ref{fact:Vopenka}, $R\in H[G_x]=\HOD^{V[G]}_x$
is a set-generic extension of $H$.
Let $\xi\in\OR$ with $\Vop^{V[G]}\sub\xi$.
So $(\xi^{+n+1})^{H[G_x]}=(\xi^{+n+1})^H$.
Since $R\in\HC^{H[G_x]}$, in $H[G_x]$ we can
iterate $R$,
and $R_{\xi^{+H}}$ has cardinality $\xi^{+H}$ there,
and $\xi^{+H}=\kappa_{\xi^{+H}}$.
So $\pow(\xi^{+H})\inter J$ also has cardinality $\xi^{+H}$ there.
So $H\not\sub J$.
\end{proof}

\subsection{No Berkeley cardinal in $\HOD_A(X)$}\label{sec:Berkeley}

We now adapt the proof of the previous section
to show that if there is a Berkeley cardinal 
(see \cite{woodin_koellner_bagaria}) then
there is no class $A$ and set $X$ such that $V=\HOD_A(X)$.
However, as mentioned in Remark \ref{rem:goldberg_usuba},
one can prove much more.

\begin{tm}\label{tm:V=HOD_A(X)_no_B}
 Let $\Omega\in\OR$ be such that $V_\Omega\sats\ZF+$``There is a Berkeley 
cardinal''. Then there is no $A\sub V_\Omega$
such that $(V_\Omega,A)\sats\ZF(A)+$``$V=\HOD_A(X)$ for some set $X$''.
\end{tm}
\begin{proof}[Proof sketch]
 Suppose otherwise and work in $(V_\Omega,A)$.
 Let $\beta$ be Berkeley.
 Let $\lambda\in\OR$, $\lambda>\beta$, be such that $(V_\lambda,A\inter 
V_\lambda)\elem_{10}(V,A)$
 and $\cof(\lambda)>\om$ (note $\beta$ is a limit of inaccessibles,
 so there are certainly ordinals of cofinality $>\om$).
 Let
 \[ j:(V_\lambda,A\inter V_\lambda)\to(V_\lambda,A\inter V_\lambda)\]
 be elementary with $\crit(j)<\beta$. By elementarity,
 there is $X\in V_\lambda$ such that $(V_\lambda,A\inter 
V_\lambda)\sats$``$V=\HOD_A(X)$'', where here $\HOD_A(X)$ is defined
locally from $A$; that is, for every $x\in V_\lambda$,
there is $\xi<V_\lambda$ such that $x$ is definable over $(V_\xi,A\inter V_\xi)$
from parameters in $X\cup\xi$.

Now we adapt the proof in the Reinhardt cardinal context
to $(V_\lambda,A\inter V_\lambda,j)$.
We have $j\rest V_\alpha\in V_\lambda$ for each $\alpha<\lambda$,
and $j``V_\lambda$ is $\in$-cofinal in $V_\lambda$.
Consider the structure $M_0=(V_\lambda,A\inter V_\lambda,j_0)$
where $j=j_0$. Define $M_1=j_0^+(M_0)$ as usual,
and note $M_1=(V_\lambda,A\inter V_\lambda,j_1)$,
where $j_1=j_0^+(j_0)$.
There are unboundedly many $\alpha<\lambda$
of cofinality $\om$ which are fixed by $j$. Fix such an $\alpha$.
Note that $k=j_0\rest V_\alpha$ is fully elementary as a map
\[ k:(V_\alpha,A\inter V_\alpha,j_0\rest V_\alpha)\to(V_\alpha,A\inter 
V_\alpha,j_0(j_0\rest V_\alpha)).\]
Therefore
$j:M_0\to M_1$
is cofinal $\Sigma_1$-elementary.

Also since these bounded fragments of $j_0$ are fully 
elementary, these facts are lifted by $j_0$, and so iterating,
$j_1:M_1\to M_2$
is also cofinal $\Sigma_1$-elementary, etc.
(In  particular,  $j_1$ also coheres $A$.)
This proceeds through all $n<\om$.
Let $i_{mn}:M_m\to M_n$ be the iteration map.

Given $\alpha<\lambda$, let $\alpha'=\sup_{n<\om}i_{0n}(\alpha)$.
Then $\alpha'<\lambda$, so note that the sequence
 $\vec{j}_\alpha=\left<i_{0n}(j\rest V_\alpha)\right>_{n<\om}\in V_\lambda$.
The map $\alpha\mapsto\vec{j}_\alpha$ is $\Sigma_1$-definable over $M_0$. 
It is now easy to adapt the usual proof to see that $M_\om$ and the maps 
$i_{n\om}:M_n\to M_\om$ are uniformly 
$\Sigma_1$-definable over $M_n$,
and $M_\om$ is wellfounded, with $\OR^{M_\om}=\lambda$.

We can now run the rest of the Reinhardt proof for a contradiction
(of course using $\HOD_A(X)$ in place of $\HOD(X)$,
which works because each $j_n$ coheres $A$).
\end{proof}

\section{Reinhardts and small forcing}\label{sec:small_forcing}

In this section we consider a different question regarding Reinhardt 
cardinals:
\begin{ques}Suppose $G$ is $(V,\PP)$-generic for some $\PP$
 in $V$, and $(V[G],j)\sats\ZFR$.
 Is there a $k$ such that $(V,k)\sats\ZFR$?
Does $j\rest V$ work?  What if we assume that $\PP\in V_{\crit(j)}$?
Note that the iterates $N_\alpha$ don't seem to help
here, since $V$ is not a proper set-generic extension of $N_\alpha$.
\end{ques}

We don't know the answer to this question in general, but we establish
some related facts here. The forcing 
arguments in this section
 relate to Laver \cite{laver_vlc}.

\begin{dfn}For a set $X$, \emph{$X$-$\AC$} is the statement
 that for every function $f$ with domain $X$,
 there is a choice function for $f$.\end{dfn}

 \begin{rem}
  The following is a straightforward adaptation of
  the standard ground definability arguments
  of Laver, Woodin and Hamkins.
Gitman and Johnstone also proved a related
fact in \cite[Main Theorem 1]{ground_def_gitman_johnstone}.
They show in particular that ($*$) if
$V\sats\ZF+\delta$-$\DC$, where $\delta$ is a cardinal,
$\card(\PP)<\delta$,
and $G$ is $(V,\PP)$-generic,
then $V$ is definable over $V[G]$ from the parameter
$\pow(\delta)^V$. Actually, they prove something more general,
in that they only assume that $\PP$ admits a gap at $\delta$,
not actually that $\card(\PP)<\delta$. We will
prove below, however, a different generalization of ($*$).
The argument is very similar
to both of those mentioned above.
 \end{rem}

 \begin{tm}\label{tm:enough_AC_ground_def}
 Let $M\sats\ZF$ and $\PP\in M$. Let $G$ be $(M,\PP)$-generic,
 and $C$ be the set of conditions of $\PP$.
 Suppose  $C$ is wellorderable in 
$M[G]$
and $M[G]\sats$``$\gamma^+$-$\DC$'' where
$\gamma=\card^{M[G]}(C)$.
 Then $M$ is definable over $M[G]$ from  the parameter
  $\pow(C^2)^M$.\footnote{The author does not
  know whether $\pow(C)^M$ suffices instead of $\pow(C^2)^M$.}
\end{tm}

This improves on ($*$) above, in that we do not assume
that $M\sats\gamma^+$-$\DC$.
(But of course,
\cite{ground_def_gitman_johnstone} only demands that $\PP$ admits
a gap at $\delta$, not that $\card(\PP)<\delta$.)

\begin{proof}
Work in $M[G]$. Let $\lambda>\eta$ be a cardinal
such that there is no surjection
$f:V_\alpha\to\lambda$ with $\alpha<\lambda$.
The following claim completes the proof:

\begin{clm} $V_\lambda\inter M$ is the unique
transitive  $N\sub V_\lambda$ such that
 (i) $\lambda=\OR\inter N$,
 (ii) $V_\alpha\inter N\in N$ for each $\alpha<\lambda$,
 (iii) $N\sats$ Separation, Pairing, Union, Powerset,
 (iv) $C\in N$ and $\pow(C^4)^N=\pow(C^4)^M$,
 (v) There is a partial order $\QQ\in\pow(C^2)^N$
 and an $(N,\QQ)$-generic $H$ such that $V_\lambda=N[H]$.
\end{clm}
\begin{proof}
We show by induction on $\xi$ that $V_\xi^N=V_\xi^M$.
This is trivial  by induction if $\xi$ is a 
limit. So assume it holds at $\xi\in[\om,\lambda)$.
Given $A\sub N$ and $B\in N$, we say that $A$ is \emph{$B$-amenable to $N$}
iff  $A\inter\rg(\pi)\in N$
for every function $\pi\in N$ such that $\pi:B\to N$.

Let $A\sub V_{\xi}$. We must show that $A\in M$ iff $A\in N$.
We do this by proving:
\begin{enumerate}
\item\label{item:A_in_M} $A\in M$ iff $A$ is $C^2$-amenable to $M$,
\item\label{item:A_in_N} $A\in N$ iff $A$ is $C^2$-amenable to $N$,
\item\label{item:pi_in_M_N} for all $\pi:C^2\to V_\xi^M=V_\xi^N$,
we have $\pi\in M$ iff $\pi\in N$.
\end{enumerate}
Putting these things together, we get $A\in M$ iff $A\in N$, as desired.

Note first that we have the Forcing Theorem
with respect to $\Sigma_0$ formulas,
and definability of the Forcing Relation for $\Sigma_0$ formulas,
over $V_\lambda^M$ and $N$; this follows from Separation
and that $\lambda$ is  closed enough.

Parts \ref{item:A_in_M},\ref{item:A_in_N}: It suffices
to consider part \ref{item:A_in_N}. For the non-trivial direction, suppose that 
$A\notin N$; we show $A$ is not $C^2$-amenable to $N$. We 
have $V_\lambda=N[H]$.
So let $\dot{A}\in N$ be such that 
$\dot{A}_H=A$. Since $A\sub V_\xi^N$ and $A\notin N$, there is $p_0\in H$
such that $N\sats p_0\forces_{\QQ}$``$\dot{A}\sub\check{V_\xi}\text{ but 
}\dot{A}\notin\check{V}$''.

Now in $V$, where $\QQ\sub C$ is wellordered and $|C|^+$-$\AC$ holds,
let $f:C\to V_\xi^N$ be such that for each $p\leq^{\QQ} p_0$, $f(p)$
is some $x$ such that
\[ N\sats p\text{ does not decide 
(w.r.t.~}\QQ)\text{ whether }\check{x}\in\dot{A}.\]

Let $\dot{f}\in N$ be such that $\dot{f}_H=f$.
Let $p_1\in H$, $p_1\leq p_0$, such that $N\sats 
p_1\forces_{\QQ}$``$\dot{f}:\check{C}\to\check{V_\xi}$,
and for each $p\leq^{\check{\QQ}}\check{p_0}$,
$\check{V}$ thinks that $p$ does not decide whether 
$\dot{f}(p)\in\dot{A}$''.

In $N$, define 
$f_*:C^2\to N$, where
for $(p,q)\in C^2$,
if $q\leq^{\QQ} p\leq^{\QQ} p_1$ then $f_*(p,q)$ is the unique $x$
such that $q\forces\dot{f}(\check{p})=\check{x}$,
if there is such an $x$ (and say $(p,q)$ is \emph{good}), and 
$f_*(p,q)=\emptyset$ otherwise.

Let $R=\rg(f_*)\in N$. We claim that $R\inter A\notin N$,
which completes the proof of part \ref{item:A_in_N}.
So let $S\in\pow(R)^N$. For each $p\leq p_1$
there is $q\leq p$ such that $(p,q)$ is good,
and so $x=f_*(p,q)\in R$; but then there are also $p'\leq p$
and $p''\leq p$ such that $p'\forces\check{x}\in\dot{A}$
and $p''\forces\check{x}\notin\dot{A}$. So by density,
there will be some such $x$ such that $x\in A$ iff $x\notin S$.
So $R\inter A\neq S$, as desired.

Part \ref{item:pi_in_M_N}:
Let $\pi^N:C^2\to V_\xi^M=V_\xi^N$, with $\pi^N\in N$.
Let $\pi^M:C^2\to V_\xi^M$, with $\pi^M\in M$.
We show that $\pi^N\in M$ and $\pi^M\in N$.

Work in $M[G]$, where $|C|=|C^2|=\gamma$ and $\gamma^+$-$\DC$ holds.
We construct a continuous increasing sequence 
$\left<X_\alpha\right>_{\alpha<\gamma^+}$ of sets
such that
$X_0=\rg(\pi^N)$,
  $X_0\cup\rg(\pi^M)\sub X_1$,
 $X_{\om\alpha+2n+1}\in M$,
  $X_{\om\alpha+2n+2}\in N$,
 $X_\alpha\sub V_\xi^M=V_\xi^N$,
 $|X_\alpha|\leq\gamma$ and $X_\alpha$ is extensional
 for each limit $\alpha$.
 
Given $X=X_{\om\alpha+2n}$,
first let $Y\sub V_{\xi}^M$ be some set with $X\sub Y$
and $|Y|=\gamma$
and such that for all $x,y\in X$, if $x\neq y$
then there is $z\in Y$ such that $z\in x$ iff $z\notin y$.
Then let $X_{\om\alpha+2n+1}\in M$
be some $Y'$ such that $Y\sub Y'\sub V_\xi^M$
and $Y'$ is the surjective image of $W$ in $M$;
this exists because $V_\lambda=V_\lambda^M[G]$,
 so $M$ can find a small covering set.
For $X_{\om\alpha+2n+2}$ it is likewise.
Make the sequence continuous.
By $\gamma^+$-$\DC$, we can proceed through $\gamma^+$ sets;
note that $\gamma$-$\DC$ ensures here that at limit stages $\alpha$,
we still have $|X_\alpha|=\gamma$.

Now since $\lambda$ is large, the sequence is in $V_\lambda$.
So let $\dot{X}^N\in N$ and $\dot{X}^M\in V_\lambda^M$
be names for it. Then we can find $p\in G$
such that for cofinally many $\alpha<\gamma^+$,
in $M$, $p$ decides the value of $\dot{X}_\alpha$,
and forces that the sequence is continuous.
Let $D^M\in M$ be the cofinal set on which $p$ decides this.
Then $\left<X_\alpha\right>_{\alpha\in D}\in M$.
Note that $D^M$ is closed.
We get a similar club $D^N\in N$.
Let $\alpha\in D^M\inter D^N$ be a limit ordinal.
Then $X_\alpha\in M\inter N$, $X_\alpha$ is extensional,
and $M,N$ both have surjections $C^2\to X_\alpha$
(since $M[G]$ has a surjection $C\to X_{\alpha+1}$,
$M$ has one $C^2\to X_{\alpha+1}$, and likewise $N$ has one 
$C^2\to X_{\alpha+2}$).

Since $X=X_\alpha$ is extensional,
we can let $\bar{X}$ be the transitive 
collapse of $X$, and $\pi:\bar{X}\to X$
the uncollapse map. Then $\bar{X},\pi\in M\inter N$
(note the recursion defining the collapse maps take less than $\lambda$ many
steps, by choice of $\lambda$).
Note that $N,M$ have surjections from $C^2$ to $\bar{X}$.
Since $\bar{X}$ is transitive, it follows that $\bar{X}$
is coded by subsets of $C^4$ in $M,N$. But $\pow(C^4)^M=\pow(C^4)^N$.
So let 
$\sigma:C^2\to\bar{X}$ be a common surjection.
Then $\pi\com\sigma:C^2\to X$ is a common surjection.
Now $\pi^M\in M$,
so $(\pi\com\sigma)^{-1}\com\pi^M\in\pow(C^4)^M=\pow(C^4)^N$,
so $\pi^M\in N$. Similarly $\pi^N\in M$.
\end{proof}
This  completes the proof of the theorem.
\end{proof}

\begin{tm}\label{tm:Reinhardt_with_enough_AC_restricts}
Let  $\eta\in\OR$ and $\PP\sub V_\eta$ and $G$ be $(V,\PP)$-generic.
Suppose $(V[G],j)\sats\ZFR$ where $\eta<\kappa=\crit(j)$,
 and $V[G]\sats$``$X$-$\AC$'' where $X=V_\eta$.
 Suppose  $j\rest V:V\to V$.
 Then $j$ is amenable to $V$.
 \end{tm}
\begin{proof}
 Suppose not and let $\alpha\in\OR$ be such that
 $j\rest V_\alpha\notin V$. Let $\beta=j(\alpha)$. Let $\widetilde{k}$
 be such that $\widetilde{k}_G=j\rest V_{\alpha+5}$.
 We may assume that
$\PP$ 
forces ``$\widetilde{k}:\check{V_{\alpha+5}}\to\check{V_{\beta+5}}$
is 
elementary
and
$\widetilde{k}\rest\check{V_\alpha}\notin\check{V}$''.

 In $V[G]$, fix a function $f:\PP\to V_\alpha^V$ such that $f(p)$ 
is some $x$ such that
 $V\sats\neg\exists y\ [ p\forces_\PP\widetilde{k}(\check{x})=\check{y}]$;
 such $f$ exists by $X$-$\AC$ in $V[G]$.
 Let $\widetilde{f}\in V$ be such that $\widetilde{f}_G=f$
 and
$V\sats\PP\forces\widetilde{f}:\check{\PP}\to\check{V_\alpha}$.
 
 In $V$, define the partial function $f_*:_{\mathrm{p}}\PP\cross\PP\to 
V_\alpha$ by
 \[ f_*(p,q)=\text{ unique }x\text{ such that 
}q\forces_\PP\widetilde{f}(\check{p})=\check{x}, \]
if $q\leq p$ and there is such an $x\in V_\alpha$, and $f^*(p,q)$ is undefined 
otherwise.
Let $D=\dom(f_*)$.

Let $g_*=j(f_*)\in V$ and let $p_0\in\PP$ be such that
$V\sats p_0\forces_\PP \widetilde{k}(\check{f_*})=\check{g}_*$.
Since $\PP\in V_{\kappa}$ and by elementarity,
$g_*:D\to V_{j(\alpha)} \text{ and } j\com f_*=g_*$.
Let $p_1\leq p_0$ with
$V\sats p_1\forces_\PP\widetilde{k}\com\check{f_*}=\check{g_*}$.
Fix $q\leq p_1$ with $(p_1,q)\in D$. Let $x=f_*(p_1,q)$. So
by choice of $f,x$,
$V\sats \neg\exists y\ [p_1\forces_\PP\widetilde{k}(\check{x})=\check{y}]$.
But then letting $y=g_*(p_1,q)$ we have
$V\sats p_1\forces_\PP 
\widetilde{k}(\check{x})=\check{y}$,
contradiction.
\end{proof}

We can now deduce that one
cannot force a ``new'' 
Reinhardt cardinal $\kappa$
such that $\kappa$-$\DC$ holds in the generic extension
with small forcing:

\begin{tm}\label{tm:Reinhardt_with_enough_AC_restricts_to_ground}
 Let $\PP\in V$ and $G$ be $(V,\PP)$-generic and $C$ the set of 
$\PP$-conditions.
Suppose
$(V[G],j)\sats\ZFR+$``$C$ is wellordered and 
 $\gamma^+$-$\DC$ holds where $\card(C)=\gamma$''.
Then  $(V,j\rest V)\sats\ZFR$.
\end{tm}
\begin{proof}
By \ref{tm:enough_AC_ground_def}, $V$ is definable over $V[G]$
from $V_{\eta+1}^V$, so
$j\rest V:V\to V$. So by \ref{tm:Reinhardt_with_enough_AC_restricts}, $j$ is 
amenable 
to $V$.
And $(V,j\rest V)$ satisfies the $\ZF$ axioms
because $(V[G],j)$ does and $(V,j\rest V)$ is definable from parameters
over $V[G]$.
\end{proof}

\begin{rem}Note that the converse direction is routine:
If $(V,j)\sats\ZFR$ and $\PP\in V_\kappa$ where $\kappa=\crit(j)$
and $G$ is $(V,\PP)$-generic, then $(V[G],j^+)\sats\ZFR$ where $j^+:V[G]\to 
V[G]$ extends $j$ with $j(G)=G$. It follows under $(V,P)\sats\ZF_2+$``
there is a super 
Reinhardt
cardinal'' then (since then there is a proper class of such) $(V[G],P[|G])$  
satisfies the same, and an easy variant gives that if $(V,P)\sats$``$\OR$ is 
total Reinhardt'', so does $(V[G],P[G])$. The other direction is less 
trivial.\end{rem}

\begin{tm}\label{tm:no_new_sR_or_tR} Assume $(V,P)\sats\ZF_2$.
Let $G$ be $(V,\PP)$-generic where $\PP\in V$. Then:
\begin{enumerate}
\item\label{item:sR} If $\kappa$ is super Reinhardt in $(V[G],P[G])$
then there is $\kappa'\in\OR^V$ which is super Reinhardt
in $(V,P)$.
\item\label{item:tR} If $\PP\in V_\delta$ and 
$(V_\delta^{V[G]},V_{\delta+1}^{V[G]})$
is total 
Reinhardt, then $(V_\delta,V_{\delta+1})$ is total Reinhardt.
\end{enumerate}
\end{tm}
\begin{proof}
Part \ref{item:tR}: Let $\delta\in\OR$ with $\PP\in V_\delta$ and 
$G$ be $(V,\PP)$-generic
and suppose that $(V_\delta^{V[G]},V_{\delta+1}^{V[G]})$ is total 
Reinhardt.
In $V$, fix $A\sub V_\delta$.
In $V[G]$, let $\kappa<\delta$ be $({<\delta},(V_\delta,A))$-reflecting
with $\PP\in V_\kappa$.
Let
\[ j:(V_\delta[G],(V_\delta,A))\to(V_\delta[G],(V_\delta,A))
\]
 be elementary 
with $\crit(j)=\kappa$.
We claim that $j\rest V_\delta\in V$
and hence, $j$ witnesses what we need in $V$.
It suffices to see that for each $\alpha<\delta$, we have $j\rest V_\alpha\in 
V$, because we have a name $\dot{j}\in V$ for $j$,
and $V_\delta$ is inaccessible, which gives $j\in V$ as usual.

So fix $\alpha<\delta$ with $\alpha=j(\alpha)$ and let $\dot{k}\in V_\delta$
be such that $\dot{k}_G=j\rest V_\alpha$. We may assume $\alpha>\kappa$.
We may assume that $\PP$ forces ``$\dot{k}:V_\alpha\to V_\alpha$
is elementary''. Let $\xi\in(\alpha+\om,\delta)$
with $\dot{k}\in V_\xi$. Note that the forcing facts in the following claim are 
defined from the relevant parameters over $V_\xi$:

\setcounter{clm}{0}
\begin{clm}
 There are $Y,\sigma,\eta\in V$, with $Y\sub V_\alpha$
 and such that (i) for each $p\in\PP$,
 if there is $x\in V_\alpha$ such that $p$ does not decide
 (with respect to $\PP$) the value of $\dot{k}(\check{x})$
then there is $x\in Y$ with this property, (ii)
 $\eta<\kappa$, and (iii) $\sigma:V_\eta\to Y$ is a surjection.
\end{clm}

\begin{proof}
(The first  part here is motivated by
Usuba's discussion of L\"owenheim-Skolem cardinals
in \cite{usuba_ls}.)
In $V[G]$, there is some $X\elem V_\xi$\footnote{We write $V_\xi$
for rank segments of $V$ and
$V_\xi^{V[G]}$ for those of $V[G]$.}
with $\PP,\dot{k}\in X$ and such that the transitive collapse
$\bar{X}$ of $X$ is in $V_\kappa^{V[G]}$. For in $V[G]$,
let
$\ell:V_\delta^{V[G]}\to V_\delta^{V[G]}$ be elementary with 
$\crit(\ell)=\kappa$ and $\ell(\kappa)>\xi$ and $\ell(V_\delta)=V_\delta$.
Then
$\ell``V_\xi\elem\ell(V_\xi)=V_{\ell(\xi)}$,
and $\ell``V_\xi\in V_\delta[G]$, and $\ell``V_\xi$ has transitive
collapse $V_\xi\in V_{j(\kappa)}^{V[G]}$. So $V_\delta[G]\sats$``There is
$X\elem V_{\ell(\xi)}$ with $\ell(\PP,\dot{k})\in X$
and such that the transitive collapse $\bar{X}$ of $X$
is in $V_{\ell(\kappa)}^{V[G]}$''. So pulling back,
in $V[G]$ we get an $X\elem V_{\xi}$ with $\PP,\dot{k}\in X$
and such that the transitive collapse $\bar{X}$ of $X$
is in $V_\kappa^{V[G]}$.
Let 
$\dot{\bar{X}}\in V$
be a name for $\bar{X}$; we may take $\dot{\bar{X}}\in V_\kappa$.
Let $\pi:\bar{X}\to X$ be the uncollapse map.
Let $\dot{\pi}\in V$ be a name for $\pi$.

Now we may assume that $\PP$ forces 
``$\dot{\pi}:\dot{\bar{X}}\to\check{V_\xi}$''.
Using the forcing relation in $V$,
looking at all outputs forceable for $\dot{\pi}$,
we can find some $Y\sub V_\alpha$ such that
$X\inter V_\alpha\sub Y$ and $\eta<\kappa$ and a surjection $\sigma:V_\eta\to Y$
which works.
\end{proof}

Let $\tau=j(\sigma)\in V$, so $\tau=j\com\sigma$. Let $p\in G$ force 
``$\dot{j}:\check{V_\delta}\to\check{V_\delta}$ is 
elementary with $\crit(\dot{j})=\check{\kappa}$
and $\dot{j}\rest\check{V_\alpha}=\dot{k}$
and $\check{\tau}=\dot{j}(\check{\sigma})=\dot{j}\com\check{\sigma}$''.
Note that $p$ forces the value of 
$\dot{k}(\check{x})=\dot{j}(\check{x})$ for every $x\in\rg(\sigma)$.
But then by the claim, $p$ forces the value of $\dot{k}(\check{x})$
for every $x\in V_\alpha$. Therefore $j\rest V_\alpha\in V$, as desired.

Part \ref{item:sR}: Let $V[G]$ be a set-generic extension such that
$\kappa$ is super-Reinhardt in $V[G]$.
Note then that if $j:V[G]\to V[G]$ is elementary with $\crit(j)=\kappa$,
then $j(\kappa)$ is also super-Reinhardt in $V[G]$
(consider $j(k)$ for various $k$).
And super Reinhardt easily implies supercompact.
So (see \cite{usuba_ls}) $V[G]\sats$``There is a proper class of supercompact 
cardinals'',
so $V[G]\sats$``There is a proper class of L\"owenheim-Skolem cardinals'',
so all set-grounds of $V[G]$ are definable from parameters over $V[G]$.
Taking $\kappa'$ super-Reinhardt then above the defining parameter,
and $j':V[G]\to V[G]$ with $\crit(j')=\kappa'$,
it follows that $j'\rest V:V\to V$. But then essentially the same
argument as above shows that $j'$ is amenable to $V$,
as required.
\end{proof}

\begin{rem}
Goldberg has shown \cite[Theorem 6.12]{goldberg_even_ordinals} that assuming
$\ZF+\DC_{V_{\lambda+1}}$ and $j:V_{\lambda+3}\to V_{\lambda+3}$ is 
$\Sigma_1$-elementary, one can prove $\Con(\ZFC+I_0)$.
It is not known at present, however, whether $\ZFR$ proves $\Con(\ZFC+I_0)$.

Assume $\ZF+k:V\to M$ is elementary with $M$ transitive, $\kappa=\crit(k)$, 
$\kappa_1=k(\kappa)$,
and $k\rest V_{\kappa_1}\in M$, and
\[ j:V_{\lambda+3}\to 
V_{\lambda+3} \text{ is }\Sigma_1\text{-elementary} \]
with 
$\kappa<\crit(j)<\kappa_\om(j)\leq\lambda<\kappa_1$.
Then there is $\PP\sub V_{\kappa}$ forcing 
$\DC_{V_{\lambda+1}}$, and
$j$ extends $\Sigma_1$-elementarily to
$j^+:V_{\lambda+3}^{V[G]}\to 
V_{\lambda+3}^{V[G]}$,
so we can apply Goldberg's result above.
To force $\DC_{V_{\lambda+1}}$, one can use a kind of argument used by Woodin 
to force $\DC$
from supercompacts.
That is, let $\PP$ be the finite support product of $\Coll(\om,V_\alpha)$,
over all $\alpha<\kappa$. Let $G$ be $(V,\PP)$-generic.
We claim $V[G]\sats\DC_{V_\alpha}$ for all $\alpha<\kappa_1$, which suffices.

For first note that for every $x\in V_{\kappa_1}$, there is 
$\kappa'<\kappa_1$
and a $\Sigma_1$-elementary $\pi:V_{\kappa}\to V_{\kappa'}$ with 
$x\in\rg(\pi)$
and $\pi(\crit(\pi))=\kappa$.
For $\pi'=k\rest V_{\kappa_1}:V_{\kappa_1}^M\to V_{\kappa'}^M$ is
$\Sigma_1$-elementary where $\kappa'=\sup k``\kappa_1$,
and $k(x)\in\rg(\pi')$, and $\pi'(\crit(\pi'))=\kappa_1$,
and $\pi'\in M$.
So the existence of such a map $\pi'$ pulls back under $k$,
yielding a $\pi$ as desired.

Now let  $R\in V_{\kappa_1}^{V[G]}$ be a finitely extendible relation,
and $\dot{R}\in V_{\kappa_1}$ with $\dot{R}_G=R$.
Let $\pi:V_\kappa\to V_{\kappa'}$ be $\Sigma_1$-elementary
with $\dot{R}_G\in\rg(\pi)$ and $\pi(\bar{\kappa})=\kappa$
 where $\bar{\kappa}=\crit(\pi)$.
 Let  $\pi(\bar{\PP},\bar{\dot{R}})=(\PP,\dot{R})$ and $\bar{G}=G\inter 
V_{\bar{\kappa}}$.
 Note that every dense subset of $\bar{\PP}$ in $V_\kappa$
 is a pre-dense subset of $\PP$, so $\bar{G}$ is $(V_\kappa,\bar{\PP})$-generic.
 Note that $\pi$ extends to a $\Sigma_1$-elementary $\pi^+:V_\kappa[\bar{G}]\to 
V_{\kappa'}[G]$
 with $\pi^+(\bar{G})=G$. Let $\bar{R}=\bar{\dot{R}}_{\bar{G}}$.
 So $\pi^+(\bar{R})=R$. But $\bar{R}$ is countable in $V[G]$,
 so $V[G]$ has an infinite branch $b$ through $\bar{R}$, and
 $\pi``b$ is a branch through $R$. 
\end{rem}

Recall that by \cite{woodin_koellner_bagaria}, $\ZF+$``$\delta_0$ is 
the least Berkeley cardinal'' implies the failure of $\cof(\delta_0)$-$\DC$ 
fails (in fact of $\cof(\delta_0)$-$\AC$).
If $j:V\to V$ is elementary then as mentioned in \cite{cumulative_periodicity} 
and \S\ref{sec:every_set_generic}, $\kappa_\om(j)$ is somewhat Berkeley. These 
and the
other considerations in this section lead to the  following conjecture:

\begin{conj}
$\ZFR+\DC$ has consistency strength
strictly greater than that of $\ZFR$.
\end{conj}

\section{Competing Interests}

The author has no competing interests.

\section{Acknowledgements}

Funded by the Deutsche 
Forschungsgemeinschaft (DFG, German Research
Foundation) under Germany's Excellence Strategy EXC 2044-390685587,
Mathematics M\"unster: Dynamics-Geometry-Structure.

The author thanks Gabriel Goldberg and Toshimichi Usuba for
their answers to certain questions and helpful feedback.

\bibliographystyle{plain}
\bibliography{../bibliography/bibliography}

\begin{thebibliography}{10}

\bibitem{woodin_koellner_bagaria}
Joan Bagaria, Peter Koellner, and W.~Hugh Woodin.
\newblock Large cardinals beyond choice.
\newblock {\em Bulletin of Symbolic Logic}, 25, 2019.

\bibitem{ground_def_gitman_johnstone}
Victoria Gitman and Thomas~A. Johnstone.
\newblock On ground model definability.
\newblock In {\em Infinity, Computability, and Metamathematics: Festschrift in
  honour of the 60th birthdays of Peter Koepke and Philip Welch}, Series:
  Tributes. College publications, London, GB, 2014.
\newblock arXiv:1311.6789.

\bibitem{virtual_lc}
Victoria Gitman and Ralf Schindler.
\newblock Virtual large cardinals.
\newblock {\em Annals of Pure and Applied Logic}, 169(12), 2018.

\bibitem{goldberg_even_ordinals}
Gabriel Goldberg.
\newblock Even ordinals and the {K}unen inconsistency.
\newblock ar{X}iv:2006.01084.

\bibitem{goldberg_email}
Gabriel Goldberg.
\newblock Personal communication, 2020.

\bibitem{cumulative_periodicity}
Gabriel Goldberg and Farmer Schlutzenberg.
\newblock Periodicity in the cumulative hierarchy.
\newblock ar{X}iv: 2006.01103.

\bibitem{gen_kunen_incon}
Joel~David Hamkins, Greg Kirmayer, and Norman~Lewis Perlmutter.
\newblock Generalizations of the kunen inconsistency.
\newblock {\em Annals of Pure and Applied Logic}, 163(12), 2012.

\bibitem{kanamori}
Akihiro Kanamori.
\newblock {\em The higher infinite: large cardinals in set theory from their
  beginnings}.
\newblock Springer monographs in mathematics. Springer-Verlag, second edition,
  2005.

\bibitem{kunen_no_R}
Kenneth Kunen.
\newblock Elementary embeddings and infinitary combinatorics.
\newblock {\em Journal of Symbolic Logic}, 36(3), 1971.

\bibitem{laver_vlc}
Richard Laver.
\newblock Certain very large cardinals are not created in small forcing
  extensions.
\newblock {\em Annals of Pure and Applied Logic}, 2007.

\bibitem{reinhardt_diss}
W.~N. Reinhardt.
\newblock {\em Topics in the metamathematics of set theory}.
\newblock PhD thesis, UC Berkeley, 1967.

\bibitem{reinhardt_remarks}
W.~N. Reinhardt.
\newblock Remarks on reflection principles, large cardinals, and elementary
  embeddings.
\newblock In {\em Axiomatic set theory (Proc. Sympos. Pure Math., Vol. XIII,
  Part II, Univ. California, Los Angeles, Calif., 1967)}, pages 189--205. Amer.
  Math. Soc., Providence, R. I., 1974.

\bibitem{vm2}
Grigor Sargsyan, Ralf Schindler, and Farmer Schlutzenberg.
\newblock Varsovian models {II}.
\newblock In preparation.

\bibitem{ZF_extenders}
Farmer Schlutzenberg.
\newblock Extenders under {ZF} and constructibility of rank-into-rank
  embeddings.
\newblock arXiv: 2006.10574.

\bibitem{con_lambda_plus_2}
Farmer Schlutzenberg.
\newblock On the consistency of {ZF} with an elementary embedding from
  {$V_{\lambda+2}$} into {$V_{\lambda+2}$}.
\newblock ar{X}iv: 2006.01077.

\bibitem{reinhardt_non-definability}
Farmer Schlutzenberg.
\newblock Reinhardt cardinals and non-definability (draft 1).
\newblock ar{X}iv: 2002.01215v1, notes.

\bibitem{suzuki_no_def_j}
Akira Suzuki.
\newblock No elementary embedding from {$V$} into {$V$} is definable from
  parameters.
\newblock {\em Journal of Symbolic Logic}, 64(4), 1999.

\bibitem{usuba_ls}
Toshimichi Usuba.
\newblock Choiceless {L}\"owenheim-{S}kolem property and uniform definability
  of grounds.
\newblock ar{X}iv: 1904.00895, 2019.

\bibitem{usuba_email}
Toshimichi Usuba.
\newblock Personal communication, 2020.

\bibitem{woodin_sem2}
W.~Hugh Woodin.
\newblock Suitable extender models {II}: Beyond $\om$-huge.
\newblock {\em Journal of Mathematical Logic}, 11(2), 2011.

\end{thebibliography}

\end{document}